\documentclass[11pt]{article}
\usepackage{}

\usepackage{amsmath,amsfonts,amsthm}
\usepackage{txfonts,wasysym,lscape,amssymb,mathrsfs,extarrows}
\usepackage[all,pdf]{xy}
\usepackage{epsfig}
\usepackage{float}      
\usepackage{multirow}
\usepackage{graphicx}
\usepackage{subfigure}
\usepackage{multirow}
\usepackage{subfigure}
\usepackage{rotating}
\usepackage{comment}
\usepackage{enumitem}
\usepackage[table,xcdraw]{xcolor}

\definecolor{pred}{RGB}{148,55,61}
\usepackage[colorlinks,linkcolor=pred,citecolor=pred,anchorcolor=pred]{hyperref}

\graphicspath{{Figures/}}

\newtheorem{theorem}{Theorem}[section]
\newtheorem{proposition}{Proposition}[section]
\newtheorem{lemma}{Lemma}[section]
\newtheorem{corollary}{Corollary}[section]

\newtheorem{remark}{Remark}[section]

\newtheorem{assumption}{Assumption}[section]

\title{ Halpern Acceleration of the Inexact Proximal Point Method of Rockafellar\thanks{Supported by National Key R\&D Program of China under project number 2022YFA1004000 and National Natural Science Foundation of China (Nos. 12371298, 12271095, and 121701185).}}

\author{Liwei Zhang
\footnote{Institute of Operations Research and Control Theory, School of Mathematical Sciences, Dalian University
of Technology, Dalian 116024, China.(lwzhang@mail.dlut.edu.cn)}, \,\, Fanli Zhuang
\footnote{Institute of Operations Research and Control Theory, School of Mathematical Sciences, Dalian University
of Technology, Dalian 116024, China.(flzhuang@mail.dlut.edu.cn)} \,\, and \,\, Ning Zhang
\footnote{School of Computer Science and Technology, Dongguan University of Technology, Dongguan 523808,China.
 (zhangning@dgut.edu.cn)}
}

\date{}
\begin{document}

\maketitle

\begin{abstract} This paper investigates a Halpern acceleration of the inexact proximal point method for solving maximal monotone inclusion problems in Hilbert spaces. The proposed Halpern inexact proximal point method (HiPPM) is shown to be globally convergent, and a unified framework is developed to analyze its worst-case convergence behavior. Under mild conditions on the inexactness tolerances, HiPPM achieves an $\mathcal{O}(1/k^{2})$ convergence rate in terms of the squared fixed-point residual. Moreover, under additional well-studied regularity conditions, the method attains a fast linear convergence rate. Building on this framework, we further extend the Halpern acceleration to the inexact augmented Lagrangian method for constrained convex optimization. In the spirit of Rockafellar’s classical results, the resulting accelerated inexact augmented Lagrangian method inherits the convergence rate and iteration complexity guarantees of HiPPM. Numerical experiments are provided to support the theoretical findings.
\vskip 6 true pt \noindent \textbf{Key words}: proximal point method, Halpern acceleration, inexactness, global convergence, rate of convergence.
\vskip 12 true pt \noindent \textbf{AMS subject classification}: 90C25, 90C30,  68Q25
\end{abstract}
\bigskip\noindent

\section{Introduction}
 \setcounter{equation}{0}

Let ${\cal H}$ be a real Hilbert space with inner product $\langle \cdot,\cdot \rangle$. The monotone inclusion problem is to
\begin{equation}\label{eq1:3}
    \text{find } z \in \mathcal{H} \text{ such that } 0 \in T(z),
\end{equation}
where $T: \mathcal{H} \rightrightarrows \mathcal{H}$ is a maximally monotone operator. This problem provides a unifying framework for a broad range of models in optimization and equilibrium theory. In particular, many convex optimization problems, variational inequalities, and saddle-point formulations can be recast as finding the zeros of a monotone operator. The proximal point method (PPM) is a fundamental approach for solving monotone inclusion problems of the form \eqref{eq1:3}. It not only enjoys strong global convergence guarantees for maximal monotone operators, but also serves as the theoretical foundation for various operator-splitting schemes. In particular, for convex optimization problems, both the augmented Lagrangian method (ALM) and the alternating direction method of multipliers (ADMM) can be interpreted as instances of the PPM when applied to their respective dual formulations.

The proximal mapping associated with parameter $c >0$, denoted by $P_{c}$ is defined by
$$
P_{c}(z)=(I+c T)^{-1}(z), z \in {\cal H}.
$$
It follows from \cite[Proposition 1]{Roc1976b} that $P_{c}(\cdot)$ is a nonexpansive mapping with
$$
0\in T(z) \mbox{ if and only if } z=P_{c}(z).
$$
Equivalently,
$$
0\in T(z) \mbox{ if and only if } Q_{c}(z) =0,
$$
where $Q_{c}=I-P_{c}$.
The classical proximal point method is based on solving a fixed point of $P_{c}$, which can be stated as follows
$$
\left\{
\begin{array}{l}
  \mbox{Choose }z^0 \in {\cal H}, \mbox{ a sequence of parameters } \{c_k\}.\\[4pt]
 \mbox{Compute }z^{k+1}=P_{c_k}(z^k) \quad \mbox{for }k=0,1,2,\ldots.
\end{array}
\right.
$$
The inexact proximal point method of Rockafellar \cite{Roc1976b} is of the following form
$$
\left\{
\begin{array}{l}
  \mbox{Choose }z^0 \in {\cal H}, \mbox{ a sequence of parameters } \{c_k\}.\\[4pt]
 \mbox{Compute }z^{k+1}\approx P_{c_k}(z^k) \mbox{ according to criterion (A) or (B)} \quad \mbox{for }k=0,1,2,\ldots.
\end{array}
\right.
$$
Two rules for the approximate calculation of $P_{c_k}(z^k)$ are criteria (A) and (B), described as follows
\begin{equation}\label{rulesAB}
\begin{array}{l}
(A) \quad \quad \quad \quad \quad \quad \|z^{k+1}-P_{c_k}(z^k)\|\leq
 \varepsilon_k, \quad \displaystyle \sum_{k=0}^{\infty} \varepsilon_k <+\infty;\\[6pt]
(B) \quad \quad \quad \quad \quad \quad \|z^{k+1}-P_{c_k}(z^k)\|\leq
 \delta_k \|z^{k+1}-z^k\|, \quad \displaystyle \sum_{k=0}^{\infty} \delta_k <+\infty.
\end{array}
\end{equation}

It was shown in \cite{brezis1978produits} that, for maximum monotone operator inclusion problems without any regularity assumptions, the proximal point method admits a worst-case convergence rate of $\mathcal{O}(1/k)$, measured in terms of the squared fixed-point residual. A tighter sublinear convergence rate was later established as $\mathcal{O}\left(1/k(1+(k-1)^{-1})^{k-1}\right) $ in \cite{gu2020tight} when the underlying Euclidean space has dimension greater than or equal to two. To further enhance the efficiency of the proximal point method, accelerating its worst-case performance has become a significant focus in both theoretical research and practical applications. In particular, Nesterov-type acceleration \cite{nesterov1983method,nesterov1988approach} has been shown to improve the rate to $\mathcal{O}(1/k^2)$. For example, Xu \cite{xu2021iteration} proposed an inexact accelerated augmented Lagrangian method for constrained convex programming.

In contrast to Nesterov-type acceleration, which was originally developed for convex optimization, Halpern’s fixed-point iteration \cite{Hal1967} was introduced for approximating fixed points of nonexpansive operators. Recent studies have revealed a close connection between accelerated proximal algorithms and Halpern’s iteration. In particular, Contreras and Cominetti \cite{contreras2023optimal} established a precise relationship between Kim’s accelerated proximal point method \cite{kim2021} and Halpern’s scheme, suggesting that Halpern-type iterations can serve as a viable mechanism for accelerating the proximal point method. More recently, a Halpern-type accelerated proximal point method with a positive semidefinite preconditioner was studied in \cite{sun2025accelerating}, this approach is shown to be equivalent to an accelerated preconditioned alternating direction method of multipliers and demonstrates strong numerical performance.

From a theoretical perspective, Lieder \cite{Lieder2021} analyzed Halpern’s iteration in Hilbert spaces and established an $\mathcal{O}(1/k^2)$ convergence rate. Moreover, under a cocoercivity assumption, Tran-Dinh \cite{tran2024halpern} showed an equivalence between Nesterov’s accelerated method and Halpern’s fixed-point iteration. Despite these advances, the existing connections are largely confined to exact algorithmic settings and do not directly extend to scenarios in which the proximal subproblems are solved inexactly, as is typically the case in large-scale applications. On the practical side, Halpern-type iterations have recently been employed to enhance GPU-based first-order methods, including preconditioned alternating direction method of multipliers, the Peaceman--Rachford method, and primal--dual hybrid gradient methods, for solving large-scale linear and quadratic programming problems; see, for example, \cite{applegate2021practical,lu2024restarted,lu2025cupdlpx,sun2025accelerating,chen2025hpr,chen2025relationships,feng2025dhpr}.


In this paper, we study a Halpern-type acceleration of the inexact proximal point method (HiPPM) for solving the maximal monotone inclusion problem~\eqref{eq1:3} in Hilbert spaces. The main contributions of this work can be summarized as follows:
\begin{itemize}
\item[1.] We establish the global convergence of the proposed inexact PPM with Halpern iteration and provide a unified framework for analyzing its convergence rate. Under mild conditions on the inexactness tolerance sequence, the squared fixed-point residual is shown to converge at the rate of~$\mathcal{O}(1/k^2)$.
\item[2.] As pointed out in~\cite{kim2021}, for strongly monotone operators, the standard proximal point method enjoys a linear convergence rate in terms of the fixed-point residual, whereas its Halpern-accelerated variant is not necessarily guaranteed to preserve such a rate. In this paper, we demonstrate that the HiPPM indeed attains linear convergence under some well-studied regularity conditions.
\item[3.] Building upon the accelerated inexact proximal point framework, we further extend the analysis to constrained convex optimization problems within the augmented Lagrangian setting. In analogy with Rockafellar’s seminal work~\cite{Roc1976b,Roc1976c}, the proposed accelerated inexact PPM naturally induces an accelerated inexact augmented Lagrangian method. This connection allows us to derive corresponding convergence rate and iteration complexity results for the accelerated iALM.
\end{itemize}

The remainder of this paper is organized as follows. Section~\ref{sec2:HiPPM} introduces the framework of the Halpern accelerated inexact proximal point algorithm in Hilbert spaces and establishes its global convergence. Section~\ref{sec3:rate-of-Conv} develops a general framework for analyzing the worst-case convergence rate of the proposed HiPPM via a direct algebraic proof, and further demonstrates its fast linear convergence under mild conditions. Section~\ref{sec4:acciALM} derives the convergence rate of the Halpern accelerated inexact augmented Lagrangian method by exploiting its close connection with the Halpern accelerated inexact proximal point algorithm. Section~\ref{sec5:numerical} presents numerical experiments that validate the theoretical results. Finally, Section~\ref{sec6:conclusion} concludes the paper.

\noindent\textbf{Notation.} Throughout the paper, $\mathcal{H}$ denotes a real Hilbert space equipped with inner product $\langle \cdot, \cdot \rangle$ and the induced norm $\|\cdot\|$. A set-valued mapping $T: \mathcal{H} \rightrightarrows \mathcal{H}$ is said to be \emph{monotone} if $\langle z - z',\, w - w' \rangle \geq 0$ whenever $w \in T(z),\, w' \in T(z')$. It is called \emph{maximal monotone} if, in addition, its graph ${\rm gph} T := \{ (z, w) \in \mathcal{H} \times \mathcal{H} : w \in T(z) \}$ is not properly contained in the graph of any other monotone operator $T': \mathcal{H} \rightrightarrows \mathcal{H}$. We denote by $\mathbb{I}_C$ the indicator function of a set $C$, and by ${\rm Dom}(f)$ the effective domain of a function $f$. All other notation will be defined as needed in the subsequent sections.

\section{Halpern acceleration of the inexact PPM}\label{sec2:HiPPM}

Building upon Rockafellar's inexact proximal point algorithm \cite{Roc1976c}, we introduce an additional Halpern-type acceleration step to enhance its convergence behavior. The resulting accelerated scheme, referred to as HiPPM, can be formulated as follows.

\begin{equation}\label{eq2:1}
\left\{
\begin{array}{l}
  \mbox{Choose }z^0 \in {\cal H}, \mbox{ a sequence of parameters } \{c_k\}.\\[4pt]
 \mbox{Compute }\bar z^{k}\approx P_{c_k}(z^k) \mbox{ according to (A) or (B)} \quad \mbox{for }k=0,1,2,\ldots;\\[4pt]
 \mbox{Compute } z^{k+1}=\displaystyle \frac{1}{k+2}z^0+\displaystyle \frac{k+1}{k+2}\bar z^k\quad \mbox{for }k=0,1,2,\ldots.
\end{array}
\right.
\end{equation}
Here criteria (A) and (B) are of the following forms:
\begin{equation}\label{rulesABa}
\begin{array}{l}
(A) \quad \quad \quad \quad \quad \quad \|\bar z^{k}-P_{c_k}(z^k)\|\leq
 \varepsilon_k, \quad \displaystyle \sum_{k=0}^{\infty} \varepsilon_k <+\infty;\\[6pt]
(B) \quad \quad \quad \quad \quad \quad \|\bar z^{k}-P_{c_k}(z^k)\|\leq
 \delta_k \|\bar z^{k}-z^k\|, \quad \displaystyle \sum_{k=0}^{\infty} \delta_k <+\infty.
\end{array}
\end{equation}

\begin{assumption}\label{ass2.1}
Let $T:{\cal H}\rightrightarrows {\cal H}$ be a maximal monotone set-valued mapping and $T^{-1}(0)\ne \emptyset$, where $T^{-1}(0)=\{z\in {\cal H}: 0 \in T(z)\}$ is the solution set of the inclusion problem (\ref{eq1:3}).
\end{assumption}

We recall the following result from \cite[Proposition~2]{Roc1976b}, which serves as a fundamental criterion for the existence of solutions to the inclusion problem $0 \in T(z)$.

\begin{lemma}\label{lem:RProp2}
Suppose that for some $\tilde z \in {\cal H}$ and $\rho >0$ one has
$$
\langle z-\tilde z, w\rangle\geq 0 \mbox{ for all }z,w \mbox{ with } w\in T(z),\|z-\tilde z\|\geq \rho.
$$
Then there exists at least one $z$ satisfying $0\in T(z)$. Moreover, this condition is both necessary and sufficient for the existence of a solution.
\end{lemma}

For notational convenience, for $k=0,1,\ldots$, define
$$
\eta^k=\bar z^k-P_{c_k}(z^k).
$$
Then, for the sequence $\{z^k\}$ generated by Algorithm (\ref{eq2:1}), one has that $\|\eta^k\|\leq \varepsilon_k$ when criterion (A) is adopted, and $\|\eta^k\|\leq \delta_k\|\bar z^k-z^k\|$ when criterion (B) is adopted.

\begin{lemma}\label{lem2.1}
Assume that Assumption \ref{ass2.1} holds and $\{z^k\}$ is generated by Algorithm (\ref{eq2:1}) with criterion (A). Then for any $z^*\in T^{-1}(0)$,
\begin{equation}\label{eq2:2b}
\|z^{k}-z^*\|\leq \|z^0-z^*\|+\displaystyle \sum_{{j=0}}^{k-1} \varepsilon_j,\,\,k=1,2,\ldots.
\end{equation}
\end{lemma}

\begin{proof}
The proof proceeds by induction. For $k=0$,
$$
z^1=\displaystyle \frac{1}{2}z^0+\displaystyle \frac{1}{2}(P_{c_0}(z^0)+\eta^0),
$$
which implies
$$
\begin{array}{ll}
\|z^1-z^*\|& =\left\|\displaystyle \frac{1}{2}[z^0-z^*]+\displaystyle \frac{1}{2}[P_{c_0}(z^0)-z^*]
+\displaystyle \frac{1}{2}\eta^0\right\|\\[8pt]
&\leq \displaystyle \frac{1}{2}\|z^0-z^*\|+\displaystyle \frac{1}{2}\|P_{c_0}(z^0)-z^*\|+\displaystyle \frac{1}{2}\|\eta^0\|\\[8pt]
&\leq \|z^0-z^*\|+\varepsilon_0.
\end{array}
$$
Suppose the result is valid until $k \geq 1$, then
$$
\begin{array}{ll}
\|z^{k+1}-z^*\|& =\left\|\displaystyle \frac{1}{k+2}[z^0-z^*]+\displaystyle \frac{k+1}{k+2}[P_{c_k}(z^k)-z^*]
+\displaystyle \frac{k+1}{k+2}\eta^k\right\|\\[8pt]
&\leq \displaystyle \frac{1}{k+2}\left\|z^0-z^*\right\|+\displaystyle \frac{k+1}{k+2}\left\|P_{c_k}(z^k)-z^*\right\|+\displaystyle \frac{k+1}{k+2}\left\|\eta^k\right\|\\[8pt]
&\leq \displaystyle \frac{1}{k+2}\left\|z^0-z^*\right\|+\displaystyle \frac{k+1}{k+2}\left\|z^k-z^*\right\|+\displaystyle \frac{k+1}{k+2}\left\|\eta^k\right\|\\[8pt]
&\leq \displaystyle \frac{1}{k+2}\left\|z^0-z^*\right\|+\displaystyle \frac{k+1}{k+2}\left(
\|z^0-z^*\|+\displaystyle \sum_{j=0}^{k-1} \varepsilon_j\right)+\displaystyle \frac{k+1}{k+2} \varepsilon_k\\[8pt]
& \leq \|z^0-z^*\|+\displaystyle \sum_{j=0}^{k} \varepsilon_j.
\end{array}
$$
This means that the result holds for $k+1$. The proof is completed.
\end{proof}

\begin{proposition}\label{prop2:1}
Let $\{z^k\}$ be any sequence generated by the Halpern accelerated proximal point algorithm under criterion (A) with $\{c_k\}$ bounded away from zero. Then $\{z^k\}$ is bounded if and only if $T^{-1}(0)\ne \emptyset$.
\end{proposition}

\begin{proof}
If $T^{-1}(0)\ne \emptyset$, then the boundedness of $\{z^k\}$ follows from Lemma \ref{lem2.1}. For the rest of the proof, we assume that $\{z^k\}$ is any bounded sequence satisfying (A). Let $s > 0$ be such that
$$
\|z^k\|<s \mbox{ and } \varepsilon_k <s \quad \forall k.
$$
Then $\{z^k\}$ has at least one weak cluster point $z^{\infty}$ and $\|z^{\infty}\|\leq s$. It follows from (\ref{eq2:1}) that
$$
z^{k+1}=\displaystyle \frac{1}{k+2}z^0+\displaystyle \frac{k+1}{k+2}\bar z^k,
$$
which implies
$$
\|\bar z^k\|=\left\|\displaystyle \frac{k+2}{k+1}z^{k+1}-\displaystyle \frac{1}{k+1}z^0\right\|
\leq 1.4 s +\left\|\displaystyle \frac{1}{k+1}z^0\right\|
< 1.5s+0.5s=2s,
$$
when $k>N_0$ for some large integer $N_0$. This implies that
$$
\|P_{c_k}(z^k)\|\leq \|\bar z^k\|+\varepsilon_k \leq 3 s, \hbox{ for } k>N_0.
$$
Then, the result follows by applying the same technique as in \cite[Theorem 1]{Roc1976b}. For completeness, the details are provided below. Consider the set-valued mapping $T'$ defined by
$$
T'(z)=T(z)+\partial h(z) \quad \forall z \in {\cal H},
$$
where
$$
{h(z)=\mathbb{I}_{3s \textbf{B}}(z)} \mbox{ where } \textbf{B}=\{z\in {\cal H}: \|z\|\leq 1\}.
$$
It is easy to obtain
$$
\partial h(z)=\left\{
\begin{array}{ll}
\{0\},  & \mbox{if } \|z\|< 3s,\\[4pt]
\{\lambda z:\lambda\geq 0\},  & \mbox{if } \|z\|= 3s,\\[4pt]
\emptyset,  & \mbox{if } \|z\|> 3s.
\end{array}
\right.
$$
Observe that $\partial h$ is a maximal monotone operator, because $h$ is a lower semicontinuous proper convex function; its effective domain is
$$
{\rm Dom}\, \partial h =3s \textbf{B}.
$$
Furthermore,
\begin{equation}\label{eq-R2.8}
{T'(z)=T(z) \mbox{ when } z \in {\rm int}\,3s \textbf{B}.}
\end{equation}
Since $\|P_{c_k}(z^k)\|<3 s$, while
$$
c_k^{-1}(z^k-P_{c_k}(z^k))\in T(P_{c_k}(z^k)).
$$
So we have
\begin{equation}\label{eq-R2.9}
P_{c_k}(z^k) \in {\rm Dom}\,(T)\cap {\rm int}\, {\rm Dom}\, (\partial h),\quad \forall k>N_0.
\end{equation}
\begin{equation}\label{eq-R2.10}
P_{c_k}(z^k) \in (I+c_k T')^{-1}(z^k), \quad \forall k >N_0.
\end{equation}
Inasmuch as ${\rm Dom}\,(T)\cap {\rm int}\, {\rm Dom}\, (\partial h)\ne \emptyset$ by (\ref{eq-R2.9}), we know that $T'$, as the sum of the maximal monotone operators $T$ and $\partial h$, is itself maximal monotone. Hence $P'_{c_k}=(I+c_k T')^{-1}$ is actually single-valued, and then \eqref{eq-R2.10} implies
$$
P_{c_k}(z^k)=P'_{c_k}(z^k) \quad \mbox{for all large } k >N_0.
$$
Thus the sequence $\{z^k\}$ can be regarded as arising from the Halpern accelerated proximal point algorithm for operator $T'$. As the effective domain ${\rm Dom}\, (T')$ is bounded, so that $(T')^{-1}(0)\ne \emptyset$ by Lemma \ref{lem:RProp2}. It follows from (\ref{eq-R2.8}) that $(T')^{-1}(0)\subseteq T^{-1}(0)$, this implies $T^{-1}(0)\ne \emptyset$. The proof is completed.
\end{proof}

\begin{theorem}\label{th2.1}
Assume that Assumption \ref{ass2.1} holds. Let $\{z^k\}$ and $\{\bar{z}^k\}$ be sequences generated by the Halpern accelerated proximal point algorithm under criterion (A) with $c_k$ bounded away from zero. Then $\{z^k\}$ converges in the weak topology to a point $z^{\infty}$ satisfying $0 \in T(z^{\infty})$ and
$$
\lim_{k\rightarrow \infty} \|Q_{c_k}(z^k)\|=0,\quad \lim_{k\rightarrow \infty}\|\bar z^k-z^k\|= 0.
$$
\end{theorem}

\begin{proof} It follows from Proposition \ref{prop2:1} that the sequence $\{z^k\}$ is a bounded sequence because Assumption \ref{ass2.1} holds. Let $s > 0$ be such that
$$
\|z^k\|<s \mbox{ and } \varepsilon_k <s, \quad \forall k.
$$
Just like the proof of Proposition \ref{prop2:1}, there exists an integer $N_0>0$ such that
$$
\|\bar z^k\|\leq 2s, \quad \|P_{c_k}(z^k)\|\leq 3s,\quad \mbox{ when }k >N_0.
$$
Let $z^{\infty}$ be a weak cluster point of $\{z^k\}$ and choose any $\bar z\in T^{-1}(0)$. It is easy to obtain the following inequalities
\begin{equation}\label{eq:R2.11a}
\|P_{c_k}(z^k)-z^0\|^2+\|Q_{c_k}(z^k)\|^2\leq\|z^k- z^0\|^2, \quad \forall k
\end{equation}
and
\begin{equation}\label{eq:R2.11}
\|P_{c_k}(z^k)-\bar z\|^2+\|Q_{c_k}(z^k)\|^2\leq\|z^k-\bar z\|^2, \quad \forall k.
\end{equation}
From (\ref{eq:R2.11a}), we have
$$
\begin{array}{l}
\|Q_{c_k}(z^k)\|^2-\|z^k-z^0\|^2+\|z^{k+1}-z^0\|^2\\[8pt]
\leq \|z^{k+1}-z^0\|^2-\|P_{c_k}(z^k)-z^0\|^2\\[8pt]
=\left(\displaystyle \frac{k+1}{k+2}\right)^2\|z^0-\bar z^k\|^2-\|(P_{c_k}(z^k)-\bar z^k)+(\bar z^k-z^0)\|^2\\[8pt]
\leq -\|(P_{c_k}(z^k)-\bar z^k)\|^2-2\langle (P_{c_k}(z^k)-\bar z^k),(\bar z^k-z^0)\rangle\\[8pt]
\leq 2\|P_{c_k}(z^k)-\bar z^k\|\|\bar z^k-z^0\|\leq 2(\|z^0\|+2s)\varepsilon_k.
\end{array}
$$
Hence we obtain
$$
\sum_{k=0}^{\infty}\|Q_{c_k}(z^k)\|^2 \leq 2(\|z^0\|+2s)\sum_{k=0}^{\infty}\varepsilon_k<+\infty.
$$
This implies that $Q_{c_k}(z^k)\rightarrow 0$.

Noting that
$$
\|z^k-\bar z^k\|=\|z^k-P_{c_k}(z^k)+P_{c_k}(z^k)-\bar z^k\|
\leq \|Q_{c_k}(z^k)\|+\|P_{c_k}(z^k)-\bar z^k\|,
$$
one has that $\|z^k-\bar z^k\|\rightarrow 0$. From (\ref{eq:R2.11}), we have
$$
\begin{array}{l}
\|Q_{c_k}(z^k)\|^2-\|z^k-\bar z\|^2+\|z^{k+1}-\bar z\|^2\\[8pt]
\leq \|z^{k+1}-\bar z\|^2-\|P_{c_k}(z^k)-\bar z\|^2\\[8pt]
=\left\|\displaystyle \frac{1}{k+2}(z^0-\bar z)+\displaystyle \frac{k+1}{k+2}(\bar z^k-\bar z)\right\|^2
-\|P_{c_k}(z^k)-\bar z^k+\bar z^k-\bar z\|^2\\[8pt]
\leq \displaystyle \frac{1}{(k+2)^2}\|(z^0-\bar z)\|^2
-2\langle P_{c_k}(z^k)-\bar z^k, \bar z^k-\bar z \rangle-\|P_{c_k}(z^k)-\bar z^k\|^2\\[8pt]
\leq \displaystyle \frac{1}{(k+2)^2}\|(z^0-\bar z)\|^2+2(2s+\|\bar z\|)\varepsilon_k.
\end{array}
$$
This implies that
$$
\|z^{k+1}-\bar z\|^2 \leq \|z^k-\bar z\|^2+\displaystyle \frac{1}{(k+2)^2}\|(z^0-\bar z)\|^2+2(2s+\|\bar z\|)\varepsilon_k,
$$
which because of
$$
\displaystyle \sum_{k=0}^{\infty} \varepsilon_k < +\infty
$$
implies the existence of
$$
\lim_{k\rightarrow \infty}\|z^k-\bar z\|=\mu <+\infty.
$$
Since $c_k$ is bounded away from zero, one has
\begin{equation}\label{eq:R2:14}
\lim_{k\rightarrow \infty}c_k^{-1}\|Q_{c_k}(z^k)\|=0.
\end{equation}
Observe that
\begin{equation}\label{eq:R2:15}
0\leq \langle z-P_{c_k}(z^k), w-c_k^{-1}Q_{c_k}(z^k)\rangle \mbox{ for all } w \in T(z).
\end{equation}
Noting that
$$
\begin{array}{ll}
\|z^{k+1}-P_{c_k}(z^k)\|&=\left\| \displaystyle \frac{1}{k+2}z^0+  \displaystyle \frac{k+1}{k+2}\bar z^k-P_{c_k}(z^k)\right\|\\[6pt]
&\leq \left\|\bar z^k-P_{c_k}(z^k)\right\|+\displaystyle \frac{1}{k+2}\left\|z^0-\bar z^k\right\|\\[6pt]
&\leq \varepsilon_k+\displaystyle \frac{1}{k+2}[2s+\|z^0\|],
\end{array}
$$
we have $\|z^{k+1}-P_{c_k}(z^k)\|\rightarrow 0$. Since $z^{\infty}$ is a weak cluster point of $\{z^k\}$, it is also a weak cluster point of $\{P_{c_k}(z^k)\}$. Then we get from (\ref{eq:R2:14}) and (\ref{eq:R2:15}) that
$$
0\leq \langle z-z^{\infty}, w\rangle,\quad \forall (z,w)\in {\rm gph}\,T.
$$

Arguing as in the proof of \cite[Theorem~1]{Roc1976b}, one can show that the sequence $\{z^k\}$ admits at most one weak cluster point. Hence, the entire sequence $\{z^k\}$ converges weakly to some $z^{\infty}$.
\end{proof}

\section{The rate of convergence}\label{sec3:rate-of-Conv}

Different from the commonly used technique of performance estimation problems in \cite{kim2021,contreras2023optimal}, in this section we adopt a direct algebraic proof to establish the convergence rate of the Halpern-accelerated proximal point algorithm with a fixed parameter $c_k \equiv c$. Furthermore, under the regularity condition considered in \cite{Roc1976b}, we show that the accelerated version admits a fast linear convergence rate when the sequence $\{c_k\}$ is nondecreasing.

\begin{proposition}\label{prop3.1}
Assume that Assumption \ref{ass2.1} holds and $\{z^k\}$ is any sequence generated by the Halpern accelerated proximal point algorithm under criterion (A) with $c_k\equiv c >0$. Then for $k =1,2,\ldots$
\begin{equation}\label{eq3:1}
\|z^k-P_{c}(z^k)\|\leq \sqrt{\displaystyle \frac{4\|z^0-z^*\|^2}{(k+1)^2}+\Delta_k},
\end{equation}
where
\begin{equation}\label{eq3:3}
\begin{array}{ll}
\Delta_k & =\displaystyle \frac{4k}{(k+1)^2}\left\langle z^k- z^0,
\eta^{k-1}\right \rangle
 +\displaystyle \frac{4}{(k+1)^2}\displaystyle \sum_{j=0}^{k-1}\left\langle z^0-z^{j},  \frac{j+1}{j+2}
\eta^{j}\right \rangle\\[6pt]
&\quad +\displaystyle \frac{4}{(k+1)^2}\displaystyle \sum_{j=0}^{k-1}j\left\langle z^0-z^{j},  \frac{j+1}{j+2}
\eta^{j}-\eta^{j-1}\right \rangle
+\displaystyle \frac{4k}{k+1}\left\langle z^{k}-P_{c}(z^{k}),
\eta^{k-1}\right \rangle\\[6pt]
&\quad -\displaystyle \frac{4}{(k+1)^2}\displaystyle \sum_{j=0}^{k-1}\displaystyle \frac{(j+1)^3}{j+2}\left\langle P_{c}(z^{j})-z^{j},
\eta^{j-1}-\eta^{j}\right \rangle\\[6pt]
& \quad +\displaystyle \frac{4}{(k+1)^2}\displaystyle \sum_{j=0}^{k-1}\displaystyle \frac{j+1}{j+2}\left\langle P_{c}(z^{j})-z^{j},
\eta^{j-1}\right \rangle
-\displaystyle \frac{4}{(k+1)^2}\displaystyle \sum_{j=1}^k\displaystyle \frac{j^2}{j+1}\|\eta^{j-1}\|^2.
\end{array}
\end{equation}
\end{proposition}

\begin{proof}
For $\eta^k=\bar z^k-P_{c}(z^k)$, $1\leq j \leq k$, $z^j$ is expressed as
\begin{equation}\label{eqL3}
z^j=\displaystyle \frac{1}{j+1}z^0+\displaystyle \frac{j}{j+1}(P_{c}(z^{j-1})+\eta^{j-1})
\mbox{ or } P_{c}(z^{j-1})=\displaystyle \frac{j+1}{j}z^j-\displaystyle \frac{1}{j}z^0-\eta^{j-1}.
\end{equation}
By nonexpansiveness of $P_{c}$, we have
\begin{equation}\label{eqL4}
\|P_{c}(z^k)-z^*\|^2\leq \|z^k-z^*\|^2 \quad \mbox{ for } z^* \in T^{-1}(0)
\end{equation}
and
\begin{equation}\label{eqL5}
\|P_{c}(z^{j})-P_{c}(z^{j-1})\|^2\leq \|z^{j}-z^{j-1}\|^2 \quad \mbox{ for } j=1,\ldots, k.
\end{equation}
Making a weighted sum of (\ref{eqL5}), we obtain
\begin{equation}\label{eqL6}
0\geq \displaystyle \sum_{j=1}^k j(j+1)\left(\|P_{c}(z^{j})-P_{c}(z^{j-1})\|^2-\|z^{j}-z^{j-1}\|^2 \right).
\end{equation}
In view of the second expression of (\ref{eqL3}), we have
$$
\begin{array}{l}
j(j+1)\|P_{c}(z^{j})-P_{c}(z^{j-1})\|^2\\[4pt]
=
j(j+1)\left\|z^j-P_{c}(z^j)+\displaystyle\frac{1}{j}[z^j-z^0]-\eta^{j-1}\right\|^2\\[4pt]
=j(j+1)\left[ \left\|z^j-P_{c}(z^j)+\displaystyle\frac{1}{j}[z^j-z^0]\right\|^2 +\|\eta^{j-1}\|^2-2\langle z^j-P_{c}(z^j)+\displaystyle\frac{1}{j}[z^j-z^0],\eta^{j-1}\rangle \right]\\[4pt]
=j(j+1)\|z^j-P_{c}(z^j)\|^2+2(j+1)\langle z^j-P_{c}(z^j), z^j-z^0\rangle+\displaystyle \frac{j+1}{j}\|z^j-z^0\|^2\\[8pt]
\quad -2j(j+1)\langle z^j-P_{c}(z^j)+\displaystyle\frac{1}{j}[z^j-z^0],\eta^{j-1}\rangle +j(j+1)\|\eta^{j-1}\|^2\\[6pt]
=A_1(j)+A_2(j)+A_3(j) -2j(j+1)\langle z^j-P_{c}(z^j)+\displaystyle\frac{1}{j}[z^j-z^0],\eta^{j-1}\rangle +j(j+1)\|\eta^{j-1}\|^2,
\end{array}
$$
where
$$
\begin{array}{ll}
A_1(j)=j(j+1)\|z^j-P_{c}(z^j)\|^2,\,\,
A_2(j)=2(j+1)\langle z^j-P_{c}(z^j), z^j-z^0\rangle,\,\,\,\,\\[4pt]
\hbox{and} A_3(j)=\displaystyle \frac{j+1}{j}\|z^j-z^0\|^2.
\end{array}
$$
In view of the first expression of (\ref{eqL3}), we have
$$
\begin{array}{l}
-j(j+1)\|z^j-z^{j-1}\|^2\\[4pt]
=-j(j+1)\left\|\displaystyle \frac{1}{j+1}(z^0-P_{c}(z^{j-1}))+P_{c}(z^{j-1})+\displaystyle \frac{j}{j+1}
\eta^{j-1}-z^{j-1}\right\|^2\\[8pt]
=-2j \langle z^0-P_{c}(z^{j-1}),P_{c}(z^{j-1})-z^{j-1}\rangle-j(j+1)\|P_{c}(z^{j-1})-z^{j-1}\|^2\\[8pt]
\quad -\displaystyle \frac{j}{j+1}\|z^0-P_{c}(z^{j-1})\|^2-2j(j+1)\left\langle   \displaystyle \frac{1}{j+1}(z^0-P_{c}(z^{j-1}))+P_{c}(z^{j-1})-z^{j-1}, \displaystyle \frac{j}{j+1}
\eta^{j-1} \right \rangle\\[8pt]
\quad -j(j+1)\left(\displaystyle \frac{j}{j+1}
\right)^2\|\eta^{j-1}\|^2\\[8pt]
=B_1(j)+B_2(j)+B_3(j)-j(j+1)\left(\displaystyle \frac{j}{j+1}
\right)^2\|\eta^{j-1}\|^2\\[8pt]
\quad -2j(j+1)\left\langle   \displaystyle \frac{1}{j+1}(z^0-P_{c}(z^{j-1}))+P_{c}(z^{j-1})-z^{j-1}, \displaystyle \frac{j}{j+1}
\eta^{j-1} \right \rangle,
\end{array}
$$
where
$$
\begin{array}{l}
B_1(j)=-\displaystyle \frac{j}{j+1}\|z^0-P_{c}(z^{j-1})\|^2,\,
B_2(j)=-2j \langle z^0-P_{c}(z^{j-1}),P_{c}(z^{j-1})-z^{j-1}\rangle,\\[6pt]
\hbox{and}\,\,B_3(j)=-j(j+1)\|P_{c}(z^{j-1})-z^{j-1}\|^2.
\end{array}
$$
We may write $B_1(j)$ as
 \begin{equation}\label{eqL9}
\begin{array}{ll}
B_1(j)& =-\displaystyle \frac{j}{j+1}\|z^0-P_{c}(z^{j-1})\|^2\\[6pt]
&=-\displaystyle \frac{j}{j+1}\left\|\displaystyle \frac{j+1}{j} z^0-\displaystyle \frac{j+1}{j}z^j+\eta^{j-1} \right\|^2\\[6pt]
&=-\displaystyle \frac{j+1}{j}\|z^0-z^j\|^2-2\langle z^0-z^j, \eta^{j-1}\rangle
-\displaystyle \frac{j}{j+1}\|\eta^{j-1}\|^2\\[6pt]
&=-A_3(j)-2\langle z^0-z^j, \eta^{j-1}\rangle
-\displaystyle \frac{j}{j+1}\|\eta^{j-1}\|^2.
\end{array}
\end{equation}
And we make a sum of $B_2(j)$:
$$
\displaystyle \sum_{j=1}^k B_2(j)=-\displaystyle \sum_{j=1}^k2j \langle z^0-P_{c}(z^{j-1}),P_{c}(z^{j-1})-z^{j-1}\rangle=\displaystyle \sum_{j=0}^{k-1}2(j+1) \langle z^0-P_{c}(z^{j}),z^{j}-P_{c}(z^{j})\rangle.
$$
Thus summing up $A_2(j)$ and $B_2(j)$ for $j=1,\ldots,k$ results in
\begin{equation}\label{eqL10}
\begin{array}{l}
\displaystyle \sum_{j=1}^k [A_2(j)+B_2(j)]=
2(k+1)\langle z^k-P_{c}(z^k), z^k-z^0\rangle\\[6pt]
~~~~~~~~~~~~~~~~~~~~~~~~~\quad +2\displaystyle \sum_{j=1}^{k-1}(j+1)\langle z^j-P_{c}(z^j),  z^j-P_{c}(z^j)\rangle+2\|z^0-P_{c}(z^0)\|^2.
\end{array}
\end{equation}
Shifting the index in the summation of $B_3(j)$
$$
\displaystyle \sum_{j=1}^k B_3(j)=-\displaystyle \sum_{j=1}^kj(j+1)\|P_{c}(z^{j-1})-z^{j-1}\|^2=
-\sum_{j=0}^{k-1}(j+1)(j+2)\|P_{c}(z^{j})-z^{j}\|^2
$$
and summing up $A_1(j)$ and $B_3(j)$ for $j=1,\ldots, k$ gives
\begin{equation}\label{eqL11}
\displaystyle \sum_{j=1}^k [A_1(j)+B_3(j)]=
k(k+1)\|z^k-P_{c}(z^k)\|^2-2\displaystyle \sum_{j=1}^{k-1} (j+1)\|z^j-P_{c}(z^j)\|^2-2\|z^0-P_{c}(z^0)\|^2.
\end{equation}
Thus, inserting (\ref{eqL9}), (\ref{eqL10}) and (\ref{eqL11}) in (\ref{eqL6}) leads to
\begin{equation}\label{eqL12}
\begin{array}{ll}
0\geq & k(k+1)\|z^k-P_{c}(z^k)\|^2+2(k+1)\langle z^k-P_{c}(z^k),z^k-z^0\rangle\\[6pt]
& -\displaystyle \sum_{j=1}^k2j(j+1)\left\langle z^j-P_{c}(z^j)+\displaystyle\frac{1}{j}[z^j-z^0],\eta^{j-1}\right\rangle +\displaystyle \sum_{j=1}^kj(j+1)\|\eta^{j-1}\|^2\\[6pt]
&
-2\displaystyle \sum_{j=1}^k\langle z^0-z^j, \eta^{j-1}\rangle-\displaystyle \sum_{j=1}^k\displaystyle \frac{j}{j+1}\|\eta^{j-1}\|^2-\displaystyle \sum_{j=1}^kj(j+1)\left(\displaystyle \frac{j}{j+1}
\right)^2\|\eta^{j-1}\|^2\\[8pt]
&-\displaystyle \sum_{j=1}^k2j(j+1)\left\langle   \displaystyle \frac{1}{j+1}(z^0-P_{c}(z^{j-1}))+P_{c}(z^{j-1})-z^{j-1}, \displaystyle \frac{j}{j+1}
\eta^{j-1} \right \rangle\\[12pt]
=&k(k+1)\|z^k-P_{c}(z^k)\|^2+2(k+1)\langle z^k-P_{c}(z^k),z^k-z^0\rangle\\[5pt]
& + C_1+C_2+C_3+\displaystyle \sum_{j=1}^k\displaystyle \frac{2j^2}{j+1}\|\eta^{j-1}\|^2,
\end{array}
\end{equation}
where
$$
\begin{array}{l}
C_1=-\displaystyle \sum_{j=1}^k2j(j+1)\left\langle z^j-P_{c}(z^j)+\displaystyle\frac{1}{j}[z^j-z^0],\eta^{j-1}\right\rangle ,\\[6pt]
C_2=-2\displaystyle \sum_{j=1}^k\langle z^0-z^j, \eta^{j-1}\rangle,\\[6pt]
C_3=-\displaystyle \sum_{j=1}^k2j(j+1)\left\langle   \displaystyle \frac{1}{j+1}(z^0-P_{c}(z^{j-1}))+P_{c}(z^{j-1})-z^{j-1}, \displaystyle \frac{j}{j+1}
\eta^{j-1} \right \rangle.
\end{array}
$$
We express $C_1+C_2$ as
\begin{equation}\label{eqC1C2}
\begin{array}{ll}
C_1+C_2&=-\displaystyle \sum_{j=1}^k2j(j+1)\left\langle z^j-P_{c}(z^j)+\displaystyle\frac{1}{j}[z^j-z^0],\eta^{j-1}\right\rangle+C_2\\[6pt]
&=-\displaystyle \sum_{j=1}^k2j(j+1)\left\langle z^j-P_{c}(z^j),\eta^{j-1}\right\rangle-\displaystyle \sum_{j=1}^k2(j+1)\left\langle z^j-z^0,\eta^{j-1}\right\rangle+C_2\\[6pt]
&
=D_1+D_2,
\end{array}
\end{equation}
where
$$
\begin{array}{l}
D_1=-\displaystyle \sum_{j=1}^k2j(j+1)\left\langle z^j-P_{c}(z^j),\eta^{j-1}\right\rangle\,\,\hbox{and}\,\,
D_2=-\displaystyle \sum_{j=1}^k2j\left\langle z^j-z^0,\eta^{j-1}\right\rangle.
\end{array}
$$
We further express $C_3$ as
\begin{equation}\label{eqC3}
\begin{array}{ll}
C_3& =-\displaystyle \sum_{j=1}^k2j(j+1)\left\langle   \displaystyle \frac{1}{j+1}(z^0-P_{c}(z^{j-1}))+P_{c}(z^{j-1})-z^{j-1}, \displaystyle \frac{j}{j+1}
\eta^{j-1} \right \rangle\\[6pt]
&=-\displaystyle \sum_{j=1}^k2j(j+1)\left\langle   \displaystyle \frac{1}{j+1}z^0+ \displaystyle \frac{j}{j+1}P_{c}(z^{j-1})-\displaystyle \frac{j}{j+1}z^{j-1}, \displaystyle \frac{j}{j+1}
\eta^{j-1} \right \rangle\\[6pt]
&
=E_1+E_2,
\end{array}
\end{equation}
where
$$
\begin{array}{l}
E_1=-2\displaystyle \sum_{j=1}^kj\left\langle z^0-z^{j-1},  \displaystyle \frac{j}{j+1}
\eta^{j-1}\right \rangle\,\,\hbox{and}\,\,
E_2=-2\displaystyle \sum_{j=1}^kj^2\left\langle P_{c}(z^{j-1})-z^{j-1},  \displaystyle \frac{j}{j+1}
\eta^{j-1}\right \rangle.
\end{array}
$$
For notational simplicity, we set $\eta^{-1}=0$, then $E_1$ can be expressed as
$$
\begin{array}{ll}
E_1& =-2\displaystyle \sum_{j=1}^kj\left\langle z^0-z^{j-1},  \displaystyle \frac{j}{j+1}
\eta^{j-1}\right \rangle\\[6pt]
&=-2\displaystyle \sum_{j=0}^{k-1}(j+1)\left\langle z^0-z^{j},  \displaystyle \frac{j+1}{j+2}
\eta^{j}\right \rangle\\[6pt]
&=-\displaystyle \sum_{j=0}^{k-1}2j\left\langle z^0-z^{j},
\eta^{j-1}\right \rangle
 -\displaystyle \sum_{j=0}^{k-1}2\left\langle z^0-z^{j},  \frac{j+1}{j+2}
\eta^{j}\right \rangle
-\displaystyle
\sum_{j=0}^{k-1}2j\left\langle z^0-z^{j},  \frac{j+1}{j+2}
\eta^{j}-\eta^{j-1}\right \rangle.
\end{array}
$$
Thus $D_2+E_1$ is expressed
$$
D_2+E_1 =-2k\left\langle z^k- z^0,
\eta^{k-1}\right \rangle
 -\displaystyle \sum_{j=0}^{k-1}2\left\langle z^0-z^{j},  \frac{j+1}{j+2}
\eta^{j}\right \rangle
-\displaystyle \sum_{j=0}^{k-1}2j\left\langle z^0-z^{j},  \frac{j+1}{j+2}
\eta^{j}-\eta^{j-1}\right \rangle.
$$
From the definition of $E_2$, we have
$$
\begin{array}{ll}
E_2& =-2\displaystyle \sum_{j=1}^kj^2\left\langle P_{c}(z^{j-1})-z^{j-1},  \displaystyle \frac{j}{j+1}
\eta^{j-1}\right \rangle\\[6pt]
 &=-2\displaystyle \sum_{j=1}^k\displaystyle \frac{j^3}{j+1}\left\langle P_{c}(z^{j-1})-z^{j-1},
\eta^{j-1}\right \rangle\\[6pt]
&=-2\displaystyle \sum_{j=0}^{k-1}\displaystyle \frac{(j+1)^3}{j+2}\left\langle P_{c}(z^{j})-z^{j},
\eta^{j}\right \rangle\\[6pt]
&=-2\displaystyle \sum_{j=0}^{k-1}j(j+1)\left\langle P_{c}(z^{j})-z^{j},
\eta^{j-1}\right \rangle
+2\displaystyle \sum_{j=0}^{k-1}\displaystyle \frac{(j+1)^3}{j+2}\left\langle P_{c}(z^{j})-z^{j},
\eta^{j-1}-\eta^{j}\right \rangle\\[6pt]
& \quad
+2\displaystyle \sum_{j=0}^{k-1}\left[ j(j+1)-\displaystyle \frac{(j+1)^3}{j+2}\right]\left\langle P_{c}(z^{j})-z^{j},
\eta^{j-1}\right \rangle\\[6pt]
&=-2\displaystyle \sum_{j=0}^{k-1}j(j+1)\left\langle P_{c}(z^{j})-z^{j},
\eta^{j-1}\right \rangle
+2\displaystyle \sum_{j=0}^{k-1}\displaystyle \frac{(j+1)^3}{j+2}\left\langle P_{c}(z^{j})-z^{j},
\eta^{j-1}-\eta^{j}\right \rangle\\[6pt]
& \quad -2\displaystyle \sum_{j=0}^{k-1}\displaystyle \frac{j+1}{j+2}\left\langle P_{c}(z^{j})-z^{j},
\eta^{j-1}\right \rangle.\\[6pt]
\end{array}
$$
Summing up $D_1$ and $E_2$, we obtain
\begin{equation}\label{eqD1E2}
\begin{array}{ll}
D_1+E_2 &
=-2k(k+1)\left\langle z^{k}-P_{c}(z^{k}), \eta^{k-1}\right \rangle
+2\displaystyle \sum_{j=0}^{k-1}\displaystyle \frac{(j+1)^3}{j+2}\left\langle P_{c}(z^{j})-z^{j},
\eta^{j-1}-\eta^{j}\right \rangle\\[6pt]
& \quad
-2\displaystyle \sum_{j=0}^{k-1}\displaystyle \frac{j+1}{j+2}\left\langle P_{c}(z^{j})-z^{j},
\eta^{j-1}\right \rangle.
\end{array}
\end{equation}
Combining (\ref{eqL12}) with (\ref{eqC1C2})--(\ref{eqD1E2}), we have
$$
\begin{array}{rl}
0& \geq  k(k+1)\|z^k-P_{c}(z^k)\|^2+2(k+1)\langle z^k-P_{c}(z^k),z^k-z^0\rangle
+ (C_1+C_2)+C_3+\displaystyle \sum_{j=1}^k\displaystyle \frac{2j^2}{j+1}\|\eta^{j-1}\|^2\\[6pt]
& =k(k+1)\|z^k-P_{c}(z^k)\|^2+2(k+1)\langle z^k-P_{c}(z^k),z^k-z^0\rangle
+ D_1+D_2+E_1+E_2+\displaystyle \sum_{j=1}^k\displaystyle \frac{2j^2}{j+1}\|\eta^{j-1}\|^2\\[6pt]
& =k(k+1)\|z^k-P_{c}(z^k)\|^2+2(k+1)\langle z^k-P_{c}(z^k),z^k-z^0\rangle
+ (D_1+E_2)+(D_2+E_1)+\displaystyle \sum_{j=1}^k\displaystyle \frac{2j^2}{j+1}\|\eta^{j-1}\|^2\\[6pt]
&=k(k+1)\|z^k-P_{c}(z^k)\|^2+2(k+1)\langle z^k-P_{c}(z^k),z^k-z^0\rangle-2k\left\langle z^k- z^0,
\eta^{k-1}\right \rangle\\[6pt]
&\quad
 -\displaystyle \sum_{j=0}^{k-1}2\left\langle z^0-z^{j}, \frac{j+1}{j+2}
\eta^{j}\right \rangle
-\displaystyle \sum_{j=0}^{k-1}2j\left\langle z^0-z^{j},  \frac{j+1}{j+2}
\eta^{j}-\eta^{j-1}\right \rangle-2k(k+1)\left\langle z^{k}-P_{c}(z^{k}),
\eta^{k-1}\right \rangle\\[6pt]
&\quad
+2\displaystyle \sum_{j=0}^{k-1}\displaystyle \frac{(j+1)^3}{j+2}\left\langle P_{c}(z^{j})-z^{j},
\eta^{j-1}-\eta^{j}\right \rangle
-2\displaystyle \sum_{j=0}^{k-1}\displaystyle \frac{j+1}{j+2}\left\langle P_{c}(z^{j})-z^{j},
\eta^{j-1}\right \rangle
+\displaystyle \sum_{j=1}^k\displaystyle \frac{2j^2}{j+1}\|\eta^{j-1}\|^2.
\end{array}
$$
Dividing the above inequality by $k+1$ and then adding (\ref{eqL4}), we obtain
$$
\begin{array}{ll}
0\geq&\displaystyle \frac{k+1}{2}\|z^k-P_{c}(z^k)\|^2-\displaystyle \frac{2}{k+1}\|z^0-z^*\|^2+\displaystyle \frac{2}{k+1}\left\| z^0-z^*-\displaystyle \frac{k+1}{2}(z^k-P_{c}(z^k)) \right\|^2\\[6pt]
&-2\displaystyle \frac{k}{k+1}\left\langle z^k- z^0,
\eta^{k-1}\right \rangle
 -\displaystyle \frac{2}{k+1}\displaystyle \sum_{j=0}^{k-1}\left\langle z^0-z^{j},  \frac{j+1}{j+2}
\eta^{j}\right \rangle\\[6pt]
&-\displaystyle \frac{2}{k+1}\displaystyle \sum_{j=0}^{k-1}j\left\langle z^0-z^{j},  \frac{j+1}{j+2}
\eta^{j}-\eta^{j-1}\right \rangle
-2k\left\langle z^{k}-P_{c}(z^{k}),
\eta^{k-1}\right \rangle\\[6pt]
&+\displaystyle \frac{2}{k+1}\displaystyle \sum_{j=0}^{k-1}\displaystyle \frac{(j+1)^3}{j+2}\left\langle P_{c}(z^{j})-z^{j},
\eta^{j-1}-\eta^{j}\right \rangle\\[6pt]
& -\displaystyle \frac{2}{k+1}\displaystyle \sum_{j=0}^{k-1}\displaystyle \frac{j+1}{j+2}\left\langle P_{c}(z^{j})-z^{j},
\eta^{j-1}\right \rangle
+\displaystyle \frac{2}{k+1}\displaystyle \sum_{j=1}^k\displaystyle \frac{j^2}{j+1}\|\eta^{j-1}\|^2.
\end{array}
$$
This implies
$$
\begin{array}{ll}
\|z^k-P_{c}(z^k)\|^2 &
\leq \displaystyle \frac{4}{(k+1)^2}\|z^0-z^*\|^2 
+\displaystyle \frac{4k}{(k+1)^2}\left\langle z^k- z^0,
\eta^{k-1}\right \rangle +\displaystyle \frac{4}{(k+1)^2}\displaystyle \sum_{j=0}^{k-1}\left\langle z^0-z^{j},  \frac{j+1}{j+2}
\eta^{j}\right \rangle\\[6pt]
& 
+\displaystyle \frac{4}{(k+1)^2}\displaystyle \sum_{j=0}^{k-1}j\left\langle z^0-z^{j},  \frac{j+1}{j+2}
\eta^{j}-\eta^{j-1}\right \rangle
+\displaystyle \frac{4k}{k+1}\left\langle z^{k}-P_{c}(z^{k}),
\eta^{k-1}\right \rangle\\[6pt]
&
-\displaystyle \frac{4}{(k+1)^2}\displaystyle \sum_{j=0}^{k-1}\displaystyle \frac{(j+1)^3}{j+2}\left\langle P_{c}(z^{j})-z^{j},
\eta^{j-1}-\eta^{j}\right \rangle
+\displaystyle \frac{4}{(k+1)^2}\displaystyle \sum_{j=0}^{k-1}\displaystyle \frac{j+1}{j+2}\left\langle P_{c}(z^{j})-z^{j},
\eta^{j-1}\right \rangle\\[6pt]
&
-\displaystyle \frac{4}{(k+1)^2}\displaystyle \sum_{j=1}^k\displaystyle \frac{j^2}{j+1}\|\eta^{j-1}\|^2,
\end{array}
$$
which completes the proof.
\end{proof}

We next derive an explicit upper bound for $\Delta_k$ in \eqref{eq3:1}. Let
\[
\beta_0 = \sum_{j=0}^{\infty} \varepsilon_j,
\quad
\kappa_0 = 2\|z^0 - z^*\| + 2\beta_0.
\]
From (\ref{eq2:2b}), it follows that
\begin{equation}\label{eqZk0}
\|z^k - z^0\|
\le \|z^k - z^*\| + \|z^0 - z^*\|
\le 2\|z^0 - z^*\| + \beta_0
< \kappa_0,
\end{equation}
and
\begin{equation}\label{eqPZk}
\|P_{c_k}(z^k) - z^k\|
\le 2\|z^k - z^*\|
\le 2\|z^0 - z^*\| + 2\beta_0
= \kappa_0.
\end{equation}
These inequalities provide uniform bounds that will facilitate the estimation of $\Delta_k$ in \eqref{eq3:1}.

\begin{proposition}\label{prop3.2}
Assume that Assumption \ref{ass2.1} holds and $\{z^k\}$ is any sequence generated by the Halpern accelerated proximal point algorithm under criterion (A) with $c_k\equiv c >0$. Then
\begin{equation}\label{eqDk}
\begin{array}{ll}
\Delta_k &\leq 4 \kappa_0\varepsilon_{k-1}\displaystyle \frac{k}{(k+1)^2}
+8\kappa_0\beta_0 \displaystyle \frac{1}{(k+1)^2}\\[8pt]
&\quad +8\kappa_0\displaystyle \frac{1}{(k+1)^2}\displaystyle \sum_{j=1}^{k-1}j\varepsilon_{j-1}+4\kappa_0\varepsilon_{k-1}\\[8pt]
&\quad +8\kappa_0\displaystyle \frac{1}{(k+1)^2}\displaystyle \sum_{j=1}^{k-1}\displaystyle (j+1)^2\varepsilon_{j-1}.
\end{array}
\end{equation}
\end{proposition}

\begin{proof}
Using (\ref{eqZk0}) and (\ref{eqPZk}), we can derive estimates for each term in \eqref{eq3:1}, which together yield \eqref{eqDk}. For notational simplicity, we set $\varepsilon_{-1}:=\varepsilon_0$. Specifically, the estimates can be summarized as follows:
\begin{itemize}
\item[(1)] $\displaystyle \left\langle z^k- z^0, \eta^{k-1}\right \rangle\leq \displaystyle \|z^k- z^0\|\cdot\|\eta^{k-1}\|
\leq \kappa_0\varepsilon_{k-1}$;
\item[(2)] $\displaystyle \sum_{j=0}^{k-1}\left\langle z^0-z^{j}, \frac{j+1}{j+2}\eta^{j}\right \rangle\leq \displaystyle \sum_{j=0}^{k-1}\|z^0-z^{j}\|\cdot\| \frac{j+1}{j+2}\eta^{j}\|\leq \kappa_0\beta_0$;
\item[(3)] $\displaystyle \sum_{j=0}^{k-1}j\left\langle z^0-z^{j}, \frac{j+1}{j+2}\eta^{j}-\eta^{j-1}\right \rangle \leq \displaystyle \sum_{j=0}^{k-1}j\|z^0-z^{j}\|[\|\eta^{j}\|+\|\eta^{j-1}\|]\leq 2\kappa_0\displaystyle \sum_{j=0}^{k-1}j\varepsilon_{j-1}$;
\item[(4)] $\left\langle z^{k}-P_{c}(z^{k}), \eta^{k-1}\right \rangle \leq \|z^{k}-P_{c}(z^{k})\|\cdot\| \eta^{k-1}\|\leq \displaystyle \kappa_0 \varepsilon_{k-1}$;
\item[(5)] $-\displaystyle \left\langle P_{c}(z^{j})-z^{j}, \eta^{j-1}-\eta^{j}\right \rangle \leq \| P_{c}(z^{j})-z^{j}\|[\|\eta^{j-1}\|+\|\eta^{j}\|]\leq 2\kappa_0\varepsilon_{j-1}$;
\item[(6)] $\displaystyle \sum_{j=0}^{k-1}\displaystyle \frac{j+1}{j+2}\left\langle P_{c}(z^{j})-z^{j}, \eta^{j-1}\right \rangle\leq \displaystyle \sum_{j=0}^{k-1}\displaystyle \frac{j+1}{j+2}\|P_{c}(z^{j})-z^{j}\|\| \eta^{j-1}\|\leq \displaystyle \kappa_0\displaystyle \sum_{j=0}^{k-1}\displaystyle \varepsilon_{j-1}$.

\end{itemize}
These estimates, together with the expression of $\Delta_k$ in \eqref{eq3:3}, imply the desired conclusion.
\end{proof}

To analyze the convergence rate of the proposed method under inexact computations, we next specify the tolerance sequence $\{\varepsilon_k\}$ in criterion~(A). A proper choice of $\{\varepsilon_k\}$ not only ensures the summability condition required by Theorem~\ref{th2.1}, but also allows us to obtain a quantitative convergence rate estimate. For criterion (A), we choose
\begin{equation}\label{eqeRule}
\varepsilon_k=\displaystyle \frac{1}{(k+2)^{1+\delta}}, \quad \delta>0, k=0,1,\ldots.
\end{equation}
Then
\begin{equation}\label{eqBeta}
\beta_0=\beta_0(\delta)=\displaystyle \sum_{k=0}^{\infty} \frac{1}{(k+2)^{1+\delta}}, \quad \delta>0
\end{equation}
and
\begin{equation}\label{eqkapa}
\kappa_0=\kappa_0(\delta)=2(\beta_0(\delta)+\|z^0-z^*\|), \quad \delta>0.
\end{equation}
Next, we show that, with a specific choice of the tolerance, the Halpern accelerated proximal point method achieves a convergence rate that is very close to the tight rate established in \cite{Lieder2021}.

\begin{theorem}\label{Th:3.1}
Assume that Assumption \ref{ass2.1} holds and $\{z^k\}$ is any sequence generated by the Halpern accelerated proximal point algorithm under criterion (\ref{eqeRule}) with $c_k\equiv c >0$. Then for $k=1,2,\ldots$
\begin{equation}\label{eqn3:1}
\|z^k-P_{c}(z^k)\|\leq \displaystyle \frac{2\|z^0-z^*\|}{(k+1)}+\sqrt{\Theta_k},
\end{equation}
where
$$
\Theta_k=\left\{
\begin{array}{ll}
8\kappa_0\beta_0 \displaystyle \frac{1}{(k+1)^2}+
4\kappa_0 \displaystyle \frac{1}{(k+1)^{2+\delta}}
+\displaystyle \frac{4\kappa_0(3-\delta)}{1-\delta} \displaystyle \frac{1}{(k+1)^{1+\delta}}
+\displaystyle \frac{8\kappa_0}{2-\delta} \frac{1}{(k+1)^{\delta}},
& \mbox{if } 0< \delta <1,\\[6pt]
4\kappa_0(1+2\beta_0) \displaystyle \frac{1}{(k+1)^2}+
4\kappa_0 \displaystyle \frac{1}{(k+1)^{3}}+
8\kappa_0 \displaystyle \frac{\ln (k+1)}{(k+1)^2}
+8\kappa_0 \displaystyle \frac{1}{k+1},& \mbox{if } \delta =1,\\[6pt]
8\kappa_0\left(\displaystyle \frac{1}{\delta-1}+\beta_0\right) \displaystyle \frac{1}{(k+1)^2}+
4\kappa_0 \displaystyle \frac{1}{(k+1)^{2+\delta}}+
4\kappa_0 \displaystyle \frac{1}{(k+1)^{1+\delta}}
+\displaystyle \frac{8\kappa_0}{2-\delta} \displaystyle \frac{1}{(k+1)^{\delta}},& \mbox{if } 1<\delta <2,\\[6pt]
8\kappa_0\left(1+\beta_0\right) \displaystyle \frac{1}{(k+1)^2}+
4\kappa_0 \displaystyle \frac{1}{(k+1)^{4}}+
4\kappa_0 \displaystyle \frac{1}{(k+1)^{3}}
+8\kappa_0 \displaystyle \frac{\ln (k+1)}{(k+1)^2},& \mbox{if } \delta =2,\\[6pt]
8\kappa_0\left(\displaystyle \frac{1}{\delta-1}+\frac{1}{\delta-2}+\beta_0\right) \displaystyle \frac{1}{(k+1)^2}+
4\kappa_0 \displaystyle \frac{1}{(k+1)^{2+\delta}}+
4\kappa_0 \displaystyle \frac{1}{(k+1)^{1+\delta}},
& \mbox{if } \delta >2,\\[6pt]
\end{array}
\right.
$$
where $\beta_0(\delta)$ and $\kappa_0(\delta)$ are defined by (\ref{eqBeta}) and (\ref{eqkapa}), respectively.
\end{theorem}

\begin{proof}
By substituting the specific choice of $\varepsilon_k$ in \eqref{eqeRule} into the definition of $\Delta_k$ in \eqref{eqDk}, we have
\begin{equation}\label{eq:Delta-k-bound}
\begin{array}{cl}
\Delta_k & \leq 4 \kappa_0\displaystyle \frac{1}{(k+1)^{2+\delta}}+8\kappa_0\beta_0 \displaystyle \frac{1}{(k+1)^2}
    +A+4\kappa_0\displaystyle \frac{1}{(k+1)^{1+\delta}},
\end{array}
\end{equation}
where,
$$
A := 8\kappa_0\displaystyle \frac{1}{(k+1)^2}\displaystyle \sum_{j=1}^{k-1}\left(\displaystyle \frac{1}{(j+1)^{\delta}}+ \displaystyle \frac{1}{(j+1)^{\delta-1}}\right) < 8\kappa_0\displaystyle \frac{1}{(k+1)^2} \int_1^{k+1} (t^{-\delta}+t^{1-\delta})dt.
$$

To proceed, we estimate the integrals in the above expression by considering different cases for the parameter $\delta > 0$. Specifically,
\begin{itemize}
    \item[(1)] For $0<\delta<1$, we have
    $$
   \frac{1}{(k+1)^2} \int_1^{k+1} (t^{-\delta}+t^{1-\delta})dt
   <\displaystyle \frac{1}{(1-\delta)(k+1)^{1+\delta}}+\displaystyle \frac{1}{(2-\delta)(k+1)^{\delta}};
    $$
    \item[(2)] For $\delta=1$, we have
    $$
    \frac{1}{(k+1)^2} \int_1^{k+1} (t^{-\delta}+t^{1-\delta})dt= \frac{1}{(k+1)^2} \int_1^{k+1} (t^{-1}+1)dt<\displaystyle \frac{\ln (k+1)}{(k+1)^{2}}+\displaystyle \frac{1}{(k+1)};
    $$
    \item[(3)] For $1<\delta<2$, we have
    $$
    \frac{1}{(k+1)^2} \int_1^{k+1} (t^{-\delta}+t^{1-\delta})dt < \displaystyle \frac{1}{(\delta-1)(k+1)^{2}}+ \displaystyle \frac{1}{(2-\delta)(k+1)^{\delta}};
    $$
    \item[(4)] For $\delta=2$, we have
    $$
     \frac{1}{(k+1)^2} \int_1^{k+1} (t^{-\delta}+t^{1-\delta})dt < \displaystyle \frac{1}{(k+1)^{2}}+ \frac{\ln (k+1)}{(k+1)^2};
    $$
    \item[(5)] For $\delta>2$, we have
     $$
     \frac{1}{(k+1)^2} \int_1^{k+1} (t^{-\delta}+t^{1-\delta})dt < \displaystyle \frac{1}{(\delta-1)(k+1)^{2}}+ \frac{\ln (k+1)}{(\delta-2)(k+1)^2}.
    $$
\end{itemize}
These estimates, together with \eqref{eq:Delta-k-bound}, complete the proof.
\end{proof}

\begin{corollary}\label{Cor:3.1}
Assume that Assumption \ref{ass2.1} holds and $\{z^k\}$ is any sequence generated by the Halpern accelerated proximal point algorithm under criterion (\ref{eqeRule}) with $c_k\equiv c >0$. Then for $k \in \textbf{N}_+$,
$$
\|z^k-P_{c}(z^k)\|\leq \left\{
\begin{array}{ll}
{\rm O}\left(\sqrt{\displaystyle \frac{8\kappa_0}{2-\delta}} \displaystyle \frac{1}{(k+1)^{\delta/2}}\right),
& \mbox{if } 0< \delta < 2,\\[12pt]
{\rm O}\left(\sqrt{8\kappa_0 \displaystyle \frac{\ln (k+1)}{(k+1)^2}}\right),& \mbox{if } \delta =2,\\[12pt]
{\rm O}\left(\sqrt{8\kappa_0\left(\displaystyle \frac{1}{\delta-1}+\frac{1}{\delta-2}+\beta_0\right) +4\|z^0-z^*\|^2 }\displaystyle \frac{1}{k+1}\right), & \mbox{if } \delta > 2.
\end{array}
\right.
$$
\end{corollary}

\begin{remark}
Corollary~\ref{Cor:3.1} provides explicit convergence rate estimates for the Halpern accelerated proximal point method under the tolerance rule~\eqref{eqeRule}. It shows that the residual $\|z^k - P_c(z^k)\|$, which measures the violation of the monotone inclusion, converges to zero at a sublinear rate depending on the decay parameter~$\delta$. Specifically, the rate improves as $\delta$ increases: for $0<\delta<2$, the decay follows ${\cal O}(k^{-\delta/2})$; at the critical case $\delta=2$, a logarithmic factor appears; and for $\delta>2$, the rate stabilizes at the optimal order ${\cal O}(1/k)$. These results demonstrate that, with appropriately chosen tolerances, the accelerated inexact PPM can approach the theoretically tight convergence behavior established in~\cite{Lieder2021}.
\end{remark}

We shall say that $T^{-1}$ is Lipschitz continuous at $0\in {\cal H}$ with modulus $a\geq 0$ if there exists a unique solution $\bar z$ to $0\in T(z)$ and for some $\bar \varepsilon>0$ we have
\begin{equation}\label{eq:R3.1}
\|z-\bar z\|\leq a\|w\| \mbox{ whenever } z\in T^{-1}(w) \mbox{ and } \|w\|\leq \bar \varepsilon.
\end{equation}
Under the same regularity conditions as those used by Rockafellar \cite{Roc1976c}, as discussed above, we now show that the accelerated version of the proximal point method still attains a fast linear convergence rate. Based on the results in Theorem \ref{th2.1}, we can derive the following results in a manner similar to that of \cite[Theorem 2]{Roc1976b}, and thus, we provide only a sketch of the proof.

\begin{theorem}\label{Th:4.1}
Assume that Assumption \ref{ass2.1} holds. Let $\{z^k\}$ and $\bar{z}^k$ be any sequences generated by the Halpern accelerated proximal point algorithm under criterion (B) of (\ref{rulesABa}) with $c_k$ nondecreasing ($c_k \nearrow c_{\infty}\leq +\infty$). Assume that $\{z^k\}$ is bounded and $T^{-1}$ is Lipschitz continuous at $0$ with modulus $a$; let
$$
\mu_k=\displaystyle \frac{a}{(a^2+c_k^2)^{1/2}} <1.
$$
Then $\{z^k\}$ converges strongly to $\bar z$, the unique solution to $0 \in T(z)$. Moreover, there is an index $\bar k$ such that
\begin{equation}\label{eq:R3.2}
\|\bar z^{k}-\bar z\| \leq \vartheta_k \|z^k-\bar z\| \quad \mbox{ for all } k \geq \bar k,
\end{equation}
where
$$
1> \vartheta_k=\displaystyle \frac{\mu_k+\delta_k}{1-\delta_k}\geq 0 \quad \mbox{ for all } k \geq \bar k,
$$
$$
\vartheta_k\rightarrow \mu_{\infty}\quad (\mbox{where } \mu_{\infty}=0 \mbox{ if } c_{\infty}=+\infty).
$$
\end{theorem}

\begin{proof}
The sequence $\{z^k\}$, being bounded, also satisfies criterion (A) of (\ref{rulesABa}) for $\varepsilon_k=\delta_k \|\bar z^k-z^k\|$. Thus the conclusions of Theorem \ref{th2.1} hold. Using the same technique in \cite[Theorem 2]{Roc1976b}, we can obtain that, there exists a sufficient large number $\bar{k}$ such that
 \begin{equation}\label{eq:R3.7}
\|P_{c_k}(z^k) -\bar z\| \leq \mu_k \|z^k-\bar z\| \quad \mbox{ if } k \geq \bar k.
\end{equation}
 Obviously we have
 $$
 \|\bar z^k-\bar z\|\leq \|\bar z^k-P_{c_k}(z^k)\|+\|P_{c_k}(z^k)-\bar z\|,
  $$
  where under criterion (B) of (\ref{rulesABa})
  $$
  \|\bar z^k-P_{c_k}(z^k)\|\leq \delta_k\|\bar z^k-z^k\|\leq \delta_k\|\bar z^k-\bar z\|+\delta_k\|z^k-\bar z\|.
  $$
 Therefore by (\ref{eq:R3.7}),
 $$
 \|\bar z^k-\bar z\|\leq \delta_k\|\bar z^k-\bar z\|+\delta_k\|z^k-\bar z\|+\mu_k \|z^k-\bar z\|.
 $$
 This inequality produces the one in (\ref{eq:R3.2}) if $\bar k \geq \tilde k$ is chosen so that (\ref{eq:R3.2}) holds, as is possible since $1 > \mu_k \searrow 0$ and $\delta_k \rightarrow 0$.

 Next we have from
 $$
 z^{k+1}=\displaystyle \frac{1}{k+2}z^0+\displaystyle \frac{k+1}{k+2}\bar z^k
 $$
  that
  $$
  \|z^{k+1}-\bar z\|\leq \displaystyle \frac{1}{k+2}\|z^0-\bar z\|+\displaystyle \frac{k+1}{k+2}\|\bar z^k-\bar z\| \rightarrow 0,
  $$
  when $k \rightarrow \infty$, this implies the strong convergence of $\{z^k\}$ to $\bar z$. The proof is completed.
\end{proof}

\section{Accelerated inexact augmented Lagrangian method}\label{sec4:acciALM}

Building upon the accelerated inexact proximal point method developed in the previous section, we now extend our analysis to its application in constrained convex optimization within the augmented Lagrangian framework. In analogy to Rockafellar's seminal work \cite{Roc1976b,Roc1976c}, the proposed accelerated inexact proximal point method naturally gives rise to an accelerated inexact augmented Lagrangian method. This connection enables us to establish convergence rate and iteration complexity results for the accelerated iALM in a manner parallel to the analysis of the classical inexact PPM. Throughout this section, we follow the notations in \cite{Roc1976c}.

Consider the optimization problem
  \begin{equation}\label{eq4.1}
  (P) \quad \quad \quad
  \begin{array}{ll}
  \displaystyle \min_{x \in C} & f_0(x)\\[4pt]
  {\rm s.t.} & f_j(x) \leq 0, j=1,\ldots,m,
  \end{array}
  \end{equation}
where $C \subseteq \mathbb R^n$ is a nonempty closed convex set, $f_i:\mathbb R^n \rightarrow \mathbb R$ is a lower semicontinuous convex function for $i=0,1,\ldots, m$. The Lagrange dual problem of Problem (\ref{eq4.1}) is of the form
$$
  (D) \quad \quad \quad \begin{array}{ll}
    \max & g_0(y) \\[6pt]
    {\rm s.t.} & y \geq 0,
    \end{array}
$$    
    where
$$ 
    g_0(y)=\inf_{x\in C} \left\{ f_0(x)+\displaystyle \sum_{j=1}^m y_j f_j(x)\right\}.
$$

  Define the ordinary Lagrangian function $l$ for Problem (\ref{eq4.1}) as
  $$
  l(x,y)=\left\{
  \begin{array}{ll}
  f_0(x)+\displaystyle \sum_{j=1}^m y_jf_j(x) & \mbox{if } x\in C \mbox{ and } y \in \mathbb R^m_+,\\[6pt]
  -\infty & \mbox{if } x\in C \mbox{ and } y \notin \mathbb R^m_+,\\[6pt]
  +\infty & \mbox{if } x\notin C.
  \end{array}
  \right.
  $$
 The essential objective function in the Lagrange dual problem (D) is
  $$
  g(y)=\inf_{x\in \mathbb R^n} l(x,y)=
  \left\{
  \begin{array}{ll}
  g_0(y) & \mbox{if } y \in \mathbb R^m_+,\\[6pt]
  -\infty & \mbox{otherwise}.
  \end{array}
  \right.
  $$
 It follows from \cite{Roc1976c} that the augmented Lagrangian method for solving (\ref{eq4.1}) is the proximal point algorithm to the inclusion
  \begin{equation}\label{eqRGE}
    0 \in \partial g(y).
    \end{equation}
 The proximal mapping associated with $\partial g$, denoted by $P_cg$, is defined by
 \begin{equation}\label{eqProx}
 (P_cg)(y)=[I+c\partial g]^{-1}(y),
 \end{equation}
 which is the unique solution to
  \begin{equation}\label{eqPM}
    \max_u \, g(u)-\displaystyle \frac{1}{2c}\|u-y\|^2.
    \end{equation}
 Now we analyze how to obtain $(P_cg)(y)$. For convenience of expression, we use $\tilde u=(P_cg)(y)$, then $\tilde u$ is the unique solution to Problem (\ref{eqPM}). Then for
 \begin{equation}\label{equs}
 Y(x,y,c)=\Pi_{\mathbb R^m_+}[y+cF(x)],
 \end{equation}
 we have
 \begin{equation}\label{eqALsolu}
 \begin{array}{ll}
 g(\tilde u)-\displaystyle \frac{1}{2c}\|\tilde u-y\|^2&
 =\displaystyle\max_u \, \left\{g(u)-\displaystyle \frac{1}{2c}\|u-y\|^2\right\}\\[6pt]
 &=\displaystyle\max_u \inf_{x\in \mathbb R^n}\left\{l(x,u)-\displaystyle \frac{1}{2c}\|u-y\|^2\right\}\\[6pt]
 &= \displaystyle\inf_{x\in \mathbb R^n}\displaystyle\max_u \left\{l(x,u)-\displaystyle \frac{1}{2c}\|u-y\|^2\right\}\\[6pt]
 &=\displaystyle\inf_{x\in \mathbb R^n} \left\{l(x,Y(x,y,c))-\displaystyle \frac{1}{2c}\|Y(x,y,c)-y\|^2\right\}\\[6pt]
 &=\displaystyle\inf_{x\in \mathbb R^n} L(x,y,c),
 \end{array}
 \end{equation}
 where $L$ is the augmented Lagrangian function
 $$
 L(x,y,c)=f_0(x)+\displaystyle \frac{1}{2c}\left[\|\Pi_{\mathbb R^m_+}[y+cF(x)]\|^2-\|y\|^2 \right].
$$
 Let $x_c(y)$ denote the solution to
 \begin{equation}\label{eqxc}
 \displaystyle\inf_{x\in \mathbb R^n} L(x,y,c).
 \end{equation}
 Then
 \begin{equation}\label{equE}
 (P_cg)(y)=\Pi_{\mathbb R^m_+}[y+cF(x_c(y))].
 \end{equation}
 In accelerated inexact proximal point algorithm, one needs to find
 $$
 \bar y^k \approx [P_{c_k}g](y^k).
 $$
 From (\ref{equE}), we may hope that this approximate solution $\bar y^k$ can be implemented by solving the problem
 \begin{equation}\label{eqxck}
 \displaystyle\inf_{x\in \mathbb R^n} L(x,y^k,c_k).
 \end{equation}
 approximately.

This idea was originally developed by Rockafellar \cite{Roc1976c}. Corresponding to criterion (A) in the inexact proximal point framework, Rockafellar proposed the following stopping criterion:
 \begin{equation}\label{ruleC}
 \begin{array}{l}
 (C) \quad \quad \quad \quad \quad \quad \phi_k (x^{k+1})-\inf \phi_k \leq
 \varepsilon_k^2/2c_k, \quad \displaystyle \sum_{k=0}^{\infty} \varepsilon_k <+\infty,
 \end{array}
 \end{equation}
 where
 \begin{equation}\label{eqphik}
 \phi_k(x)=L(x,y^k,c_k).
 \end{equation}
  Now we are in a position to state the accelerated augmented Lagrangian method:
 \begin{equation}\label{eqAlg}
 \left\{
 \begin{array}{l}
  \mbox{Choose }y^0 \in \mathbb R^m_+, \mbox{ a sequence of parameters } \{c_j\}.\\[4pt]
 \mbox{Compute } x^{k+1}\approx \arg\min \phi_k(x) \mbox{ according to criterion (C) for } k=0,1,2,\ldots;\\[4pt]
 \mbox{Compute } \bar y^{k}=Y(x^{k+1},y^k,c_k) \mbox{ for } k=0,1,2,\ldots;\\[4pt]
 \mbox{Compute } y^{k+1}=\displaystyle \frac{1}{k+2}y^0+\displaystyle \frac{k+1}{k+2}\bar y^k\quad \mbox{for }k=0,1,2,\ldots.
\end{array}
\right.
\end{equation}

As a direct consequence of \cite[Proposition~6]{Roc1976c}, we obtain the following result, which reveals the connection between criterion~(C) for the inexact ALM and criterion~(A) for the inexact PPM.
\begin{corollary}\label{coroRuleA}
For $P_{c}g$ as in (\ref{eqProx}), $\phi_k$ as in (\ref{eqphik}), $x^{k+1}$ is generated by
$$
x^{k+1}\approx \arg\min \phi_k(x)
$$
according to criterion (C) and $\bar y^{k}=Y(x^{k+1},y^k,c_k)$, one has
$$
 \|\bar y^k-[P_{c_k}g](y^k)\|\leq \varepsilon_k,
 $$
 namely criterion (A) is satisfied for $\partial g$.
\end{corollary}

We recall the notions of an asymptotically minimizing sequence and the asymptotic infimum of Problem~\eqref{eq4.1}, as introduced by Rockafellar~\cite{Roc1976c}. A sequence $\{x^k\} \subseteq C$ is said to be asymptotically minimizing for Problem~\eqref{eq4.1} if
\begin{equation}\label{eqR4.12}
\limsup_{k \to +\infty} f_i(x^k) \le 0, \quad i = 1, \ldots, m,
\end{equation}
and if $\limsup_{k \to +\infty} f_0(x^k)$ attains the smallest possible value among all sequences in $C$ satisfying~\eqref{eqR4.12}. This smallest value is referred to as the asymptotic infimum of Problem~\eqref{eq4.1}, and is denoted by $\operatorname{asym}\inf (P)$.

\begin{theorem}\label{th4.1}
Suppose $\sup (D)>-\infty$, and let the accelerated augmented Lagrangian method (\ref{eqAlg}) be executed. If the generated sequence $\{y^k \}\subset \mathbb R^m_+$ is bounded, then $y^k \rightarrow y^{\infty}$, where $y^{\infty}$ is some optimal solution to Problem (D), and $\{x^k\}$ is asymptotically minimizing for Problem (P) with
 \begin{equation}\label{R4.13}
 f_i(x^{k+1})\leq \displaystyle \frac{1}{c_k}\left\{y^{k+1}_i-y^{k}_i+\displaystyle \frac{1}{k+1}[y^{k+1}_i-y^0_i]\right\}\rightarrow 0, \quad \mbox{ for }i=1,\ldots,m,~~~~
 \end{equation}
 \begin{equation}\label{R4.14}
 \begin{array}{l}
 f_0(x^{k+1}) -{\rm asym}\inf (P) \\[2mm]
 \leq \displaystyle \frac{1}{2c_k}\left[\varepsilon_k^2+\|y^{k}\|^2-\|y^{k+1}\|^2-\displaystyle \frac{2k+3}{(k+1)^2}\left\langle y^{k+1}-y^0, y^{k+1}-\displaystyle \frac{1}{2k+3}y^0\right\rangle\right].
 \end{array}
 \end{equation}

 The boundedness of $\{y^k\}$ under (C) is actually equivalent to the existence of an optimal solution to Problem (D). It holds if Problem (P) satisfies the Slater condition; in this case one has $\max\,(D)=\inf \,(P)={\rm asym}\, \inf (P)$.

 If $\{y^k\}$ is bounded and there exists an $\alpha_0$ such that the set
 $$
 \{x\in \mathbb R^n: x \mbox{ is feasible for } (P) \mbox{ and } f_0(x) \leq \alpha_0\}
 $$
 is nonempty and bounded, then the sequence $\{x^k\}$ is also bounded, and all of its cluster points are optimal solutions to Problem (P).
\end{theorem}
\begin{proof}
 Corollary \ref{coroRuleA} shows that criterion (C) implies criterion (A) for $T=-\partial g$. Thus the accelerated augmented Lagrangian method is reduced to the Halpern accelerated proximal point method under criterion (A) to $0\in \partial -g(y)$. Then from Theorem \ref{th2.1}, we have that the sequence $\{y^k\}$ converges to a solution $y^{\infty}$ to $0\in -\partial g(y)$, namely an optimal solution to Problem (D). From the definition of $y^{k+1}$, we have
\begin{equation}\label{eqykp1}
\bar y^k=\displaystyle \frac{k+2}{k+1}y^{k+1}-\displaystyle \frac{1}{k+1}y^0
\end{equation}
and thus
\begin{equation}\label{eqbaryk}
\bar y^k_i- y^k_i= y^{k+1}_i-y^{k}_i+\displaystyle \frac{1}{k+1}[y^{k+1}_i-y^0_i]
\end{equation}
and
$$
\bar y^{k}_i-y^k_i=Y_i(x^{k+1},y^k,c_k)=
\max \left\{-y^k_i,c_k f_i(x^{k+1})\right\}\geq c_kf_i(x^{k+1})
$$
and $y^{k+1}-y^k \rightarrow 0$ as well as the boundedness of $\{y^k\}$, we obtain (\ref{R4.13}) from (\ref{eqbaryk}).

Now we prove (\ref{R4.14}). Observing that
$$
\phi_k(x)=f_0(x)+\displaystyle \frac{1}{2c_k}[\|Y(x,y^k,c_k\|^2-\|y^k\|^2]\quad \mbox{ for } x\in C,
$$
we have
\begin{equation}\label{eqR4.16}
\phi_k(x^{k+1})-f_0(x^{k+1})=\displaystyle \frac{1}{2c_k}[\|\bar y^k\|^2-\|y^k\|^2].
\end{equation}
From the definition of $[P_{c_k}g](y^k)$, we have
$$
\begin{array}{ll}
\displaystyle \inf_x \phi_k(x)& =\inf_x L(x,y^k,c_k)\\[6pt]
& =\max_y\left\{g(y)-\displaystyle \frac{1}{2c_k}\|y-y^k\|^2\right\}\\[6pt]
&=g([P_{c_k}g](y^k))-
\displaystyle \frac{1}{2c_k}\|[P_{c_k}g](y^k)-y^k\|^2,
\end{array}
$$
which implies
\begin{equation}\label{eqR4.17}
\inf \phi_k \leq g([P_{c_k}g](y^k))\leq \max \, (D).
\end{equation}
Combining (\ref{eqR4.16}) and (\ref{eqR4.17}) we get
\begin{equation}\label{eqR4.18}
\begin{array}{ll}
f_0(x^{k+1})-\max (D) &\leq \phi_k(x^{k+1})-\inf \phi_k +\displaystyle\frac{1}{2c_k}[\|y^k\|^2-\|\bar y^{k}\|^2]\\[6pt]
&\leq \displaystyle\frac{1}{2c_k}[\varepsilon_k^2+\|y^k\|^2-\|\bar y^{k}\|^2].
\end{array}
\end{equation}
In view of (\ref{eqykp1}), we have
\begin{equation}\label{eqykp12}
\begin{array}{ll}
\|\bar y^k\|^2& =\left\|y^{k+1}+\displaystyle \frac{1}{k+1}y^{k+1}-\displaystyle \frac{1}{k+1}y^0\right\|^2\\[7pt]
&=\|y^{k+1}\|^2+\displaystyle \frac{2k+3}{(k+1)^2} \left\langle y^{k+1}-y^0, y^{k+1}-\displaystyle \frac{1}{2k+3}y^0\right\rangle.
\end{array}
\end{equation}
But every $x \in C$ satisfies
$$
f_0(x)+\displaystyle \sum_{i=1}^{\infty} y^{\infty}_i f_i(x)\geq \inf_x l(x,y^{\infty})=
g(y^{\infty})=\max \, (D),
$$
so that $\max \, (D) \leq {\rm asym}\,\inf \,(P)$. Therefore (\ref{eqR4.18}) and (\ref{eqykp12}) imply (\ref{R4.14}).
\end{proof}

Now we discuss the complexity of the accelerated augmented Lagrangian method. Define the average point
\begin{equation}\label{eqxtil}
\tilde{x}^k=\displaystyle \sum_{j=1}^kc_j x^j, \quad
\mbox{ where } c_j=c_{j-1}\left(\displaystyle \sum_{l=1}^kc_{l-1}\right)^{-1}.
\end{equation}

\begin{proposition}\label{prop4.1}
Suppose the solution set to Problem (D) is nonempty, and let the accelerated augmented Lagrangian method (\ref{eqAlg}) be executed. Then the generated sequences $\{x^k\}\subseteq C$ and $\{y^k \}\subset \mathbb R^m_+$ satisfy
 \begin{equation}\label{R4.13a}
 f_i(\tilde x^{k})\leq \displaystyle \left(\displaystyle \sum_{j=0}^{k-1}c_j\right)^{-1}\left\{y^{k}_i-y^{0}_i+\displaystyle \sum_{j=1}^k \frac{2 \Delta_0}{j}\right\}, \quad \mbox{ for }i=1,\ldots,m,
 \end{equation}
 \begin{equation}\label{R4.14a}
 f_0(\tilde x^{k}) -{\rm asym}\inf (P) \leq \displaystyle \left(\displaystyle \sum_{j=0}^{k-1}2c_j\right)^{-1}\left[\displaystyle \sum_{j=0}^{k-1}\varepsilon_j^2+\|y^{0}\|^2-\|y^{k}\|^2+\displaystyle \sum_{j=1}^{k}\displaystyle \frac{12\Delta_0^2}{j}\right],
 \end{equation}
 where $\Delta_0$ is a constant satisfying
 \begin{equation}\label{eqdelta0}
 \|y^k\|\leq \Delta_0 \quad \mbox{ for } k=0,1,2,\ldots.
 \end{equation}
 \end{proposition}
\begin{proof}
First of all, in view of Theorem \ref{th4.1}, as the solution set to Problem (D) is nonempty, the sequence $\{y^k\}$ is bounded and there exists a positive constant $\Delta_0$ satisfying (\ref{eqdelta0}). From (\ref{R4.13}) in Theorem \ref{th4.1}, we have
\begin{equation}\label{R4.1300}
 c_{k-1}f_i(x^{k})\leq \left\{y^{k}_i-y^{k-1}_i+\displaystyle \frac{1}{k}2\Delta_0\right\} \quad \mbox{ for }i=1,\ldots,m,
 \end{equation}
which implies that
$$
\displaystyle \sum_{j=1}^k c_{j-1}\times \displaystyle \sum_{j=1}^k c_jf_i(x^{j})\leq \left\{y^{k}_i-y^0_i+\displaystyle \sum_{j=1}^k\frac{1}{j}2\Delta_0\right\} \quad \mbox{ for }i=1,\ldots,m.
$$
From Jensen's inequality, we obtain (\ref{R4.13a}) from the above inequality.

 In view of (\ref{R4.14}), $(2k+3)/(k+1)^2\leq 3/(k+1)$ and
 $$
 -\left\langle y^{k+1}-y^0, y^{k+1}-\displaystyle \frac{1}{2k+3}y^0\right\rangle
 \leq 4\Delta_0^2,
 $$
 we obtain
 \begin{equation}\label{R4.14000}
 c_{j-1}[f_0(x^{j}) -{\rm asym}\inf (P) ]\leq \displaystyle \frac{1}{2}\left[\varepsilon_{j-1}^2+\|y^{j-1}\|^2-\|y^{j}\|^2-\displaystyle \frac{3}{j}\times 4\Delta_0^2\right],
 \end{equation}
 which implies
 $$
 \displaystyle \sum_{j=1}^k c_{j-1}\times\left[\displaystyle \sum_{j=1}^k c_jf_0(x^{j}) -{\rm asym}\inf (P) \right]\leq \displaystyle \frac{1}{2}\left[\displaystyle \sum_{j=1}^k\varepsilon_{j-1}^2+\|y^{0}\|^2-\|y^{k}\|^2-\displaystyle \sum_{j=1}^k \frac{12}{j}\Delta_0^2\right].
 $$
 Again from Jensen's inequality, we obtain (\ref{R4.14a}) from the above inequality.
\end{proof}

Based on Proposition~\ref{prop4.1}, by considering different choices of the sequence $\{c_k\}$, we can further establish the convergence rate and iteration complexity of the accelerated inexact augmented Lagrangian method. The corresponding results are summarized in the following theorem.

\begin{theorem}\label{th4.2}
Suppose the solution set to Problem (D) is nonempty, and let the accelerated augmented Lagrangian method (\ref{eqAlg}) be executed in which $\varepsilon_k$ is defined by (\ref{eqeRule}), and $\{x^k\}\subseteq C$ and $\{y^k \}\subset \mathbb R^m_+$ are generated sequences.
 \begin{itemize}[leftmargin=2em]
 \item[{\rm (1)}] If $c_k$ is a sequence of increasing positive parameters with $c_0>0$, then
 \begin{equation}\label{R4.13aconCom}
 f_i(\tilde x^{k})\leq \displaystyle \frac{1}{c_0k}[3\Delta_0+2\ln k]\quad \mbox{ for }i=1,\ldots,m,
 \end{equation}
 \begin{equation}\label{R4.14aobjCom}
 f_0(\tilde x^{k}) -{\rm asym}\inf (P) \leq \displaystyle \frac{1}{2c_0k}\left[\displaystyle \frac{2(1+\delta)}{1+2\delta}+13\Delta_0^2+12 \Delta_0^2 \ln k\right]
 \end{equation}
 \item[{\rm (2)}] If $c_k=c_0(k+1)$ with $c_0>0$, then
 \begin{equation}\label{R4.13aconComa}
 f_i(\tilde x^{k})\leq \displaystyle \frac{2}{c_0k(k+1)}[3\Delta_0+2\ln k]\quad \mbox{ for }i=1,\ldots,m,
 \end{equation}
 \begin{equation}\label{R4.14aobjComa}
 f_0(\widetilde x^{k}) -{\rm asym}\inf (P) \leq \displaystyle \frac{1}{c_0k(k+1)}\left[\displaystyle \frac{2(1+\delta)}{1+2\delta}+13\Delta_0^2+12 \Delta_0^2 \ln k\right]
 \end{equation}
 \end{itemize}
 where $\Delta_0$ is a constant (\ref{eqdelta0}).
 \end{theorem}

\begin{proof}
 From (\ref{R4.13a}), we have for $i=1,\ldots,m$ that
$$
 f_i(\tilde x^{k}) \leq \displaystyle \left(\displaystyle \sum_{j=0}^{k-1}c_j\right)^{-1}\left[3\Delta_0+\displaystyle \sum_{j=2}^k \frac{2 \Delta_0}{j}\right]\leq \left(k c_0\right)^{-1}\left[3\Delta_0+\displaystyle \int_1^k \frac{2 \Delta_0}{t}dt\right] \leq \displaystyle \frac{1}{c_0k}[3\Delta_0+2\ln k],
 $$
 which is just (\ref{R4.13aconCom}). From (\ref{R4.14a}) and (\ref{eqeRule}), we have that
 $$
 \begin{array}{ll}
 f_0(\tilde x^{k}) -{\rm asym}\inf (P) & \leq \displaystyle \left(\displaystyle \sum_{j=0}^{k-1}2c_j\right)^{-1}\left[\displaystyle \sum_{j=0}^{k-1}\displaystyle \frac{1}{(j+1)^{2+2\delta}}+13\Delta_0^2+\displaystyle \sum_{j=2}^{k}\displaystyle \frac{12\Delta_0^2}{j}\right]\\[16pt]
 &\leq \displaystyle \left(2c_0k \right)^{-1}\left[\displaystyle \sum_{j=1}^{k-1}\displaystyle \frac{1}{(j+1)^{2+2\delta}}+1+13\Delta_0^2+\displaystyle \sum_{j=2}^{k}\displaystyle \frac{12\Delta_0^2}{j}\right]\\[16pt]
 &\leq \displaystyle \left(2c_0k \right)^{-1}\left[\displaystyle \int_1^k\displaystyle \frac{1}{t^{2+2\delta}}dt+1+13\Delta_0^2+\displaystyle \int_1^{k}\displaystyle \frac{12\Delta_0^2}{t}dt\right]\\[16pt]
 & \leq \displaystyle \left(2c_0k \right)^{-1}\left[\displaystyle \frac{1}{1+2\delta}+1+13\Delta_0^2+\displaystyle 12\Delta_0^2\ln k\right]\\[16pt]
 &=\left(2c_0k \right)^{-1}\left[\displaystyle \frac{2(1+\delta)}{1+2\delta}+13\Delta_0^2+12 \Delta_0^2 \ln k\right].
 \end{array}
 $$
 This proves (\ref{R4.14aobjCom}). Inequalities (\ref{R4.13aconComa}) and (\ref{R4.14aobjComa}) are easily obtained in a similar way.
 \end{proof}

\section{Numerical Experiments}\label{sec5:numerical}

In this section, we report numerical experiments designed to validate the theoretical results established in this paper. In particular, we investigate the performance of the inexact proximal point method with Halpern accelerated iteration (HiPPM) on two classes of optimization problems with distinct regularity properties. The first example is the nuclear norm regularized least squares problem \cite{ma2011fixed}, which does not necessarily satisfy the regularity assumptions required for the linear convergence analysis. The second example is the sparse logistic regression problem, which satisfies metric subregularity and hence admits a fast linear convergence rate for the inexact proximal point method even in the absence of Halpern acceleration. All experiments are implemented in Matlab 2018b.

\subsection{Nuclear norm regularized Least square problem}

The following nuclear norm regularized least squares problem is commonly used as a convex surrogate for the rank minimization problem:
\begin{equation}\label{model:nuclear_ls}
\min_{z\in\mathbb{R}^{m\times n}} \; \frac{1}{2}\|Az-B\|^2+\lambda\|z\|_{*},
\end{equation}
where $\lambda>0$ and $\|z\|_{*}$ denotes the nuclear norm of $z$, defined as the sum of its singular values.

In this experiment, we adopt the stopping criterion (A) with $\varepsilon_k$ specified by \eqref{eqeRule} to validate the main result stated in Theorem~\ref{Th:3.1}. Recall that the stopping criterion is given by
\begin{equation}\label{eq:stop_criterion}
\|\bar{z}^k-P_c(z^k)\| \leq \varepsilon_k=\frac{1}{(k+2 )^{1+\delta}}.
\end{equation}
Here, $\bar{z}^k$ is an approximate solution obtained by solving the $k$ th subproblem of the HiPPM,
\begin{equation}\label{model:subproblem}
\min_{z\in\mathbb{R}^{m\times n}} \; \Phi_k(z)
:= \frac{1}{2}\|Az-B\|^2+\lambda\|z\|_{*}+\frac{1}{2c}\|z-z^k\|^2.
\end{equation}

We solve the subproblem \eqref{model:subproblem} using a dual-based accelerated gradient descent (AGD) method. Specifically, AGD is applied to the corresponding dual problem
\begin{equation}\label{model-LS-dual}
\max_{x} \; \Psi_k(x)
:= -\frac{1}{2}\|x\|^2-\langle B-Az^k, x\rangle
+c^{-1}\mathcal{M}_{c\varphi}(z^k-cA^Tx)
-\frac{c}{2}\|A^Tx\|^2,
\end{equation}
where $\mathcal{M}_{c\varphi}(\cdot)$ denotes the Moreau envelope of $\varphi(z):=\lambda\|z\|_{*}$.

Following the discussion in \cite{zhang2020proximal}, the stopping criterion involving the unknown proximal point $P_c(Z^k)$ can be replaced by the following computable condition:
\begin{equation}\label{eq:computable_stop}
\|\bar{z}^k-P_c(z^k)\|
\leq \sqrt{2c\bigl(\Phi_k(\bar{z}^k)-\Psi_k(x^{k+1})\bigr)}
\leq \varepsilon_k,
\end{equation}
where $x^{k+1}$ is an approximate solution to the dual problem \eqref{model-LS-dual}, and
\[
\bar{z}^k = cP_c\bigl(c^{-1}z^k-A^Tx^{k+1}\bigr).
\]

In the experiment, we set $m=300$ and $n=200$, and generate a ground truth matrix with rank $r=50$. The matrices $A\in\mathbb{R}^{50\times300}$ and $B\in\mathbb{R}^{50\times200}$ are generated randomly. The regularization parameter is chosen as $\lambda=1$, and the proximal parameter is fixed at $c=10$. Since the Halpern iteration is sensitive to the choice of the anchor point $Z^0$, the selection of the initial point plays a critical role in the convergence behavior. To mitigate this effect, we adopt a heuristic restarting strategy similar to that used in \cite{kim2021}. All reported results are averaged over 20 independent runs.

\begin{figure}[tbp]
\centering
\includegraphics[scale=0.33]{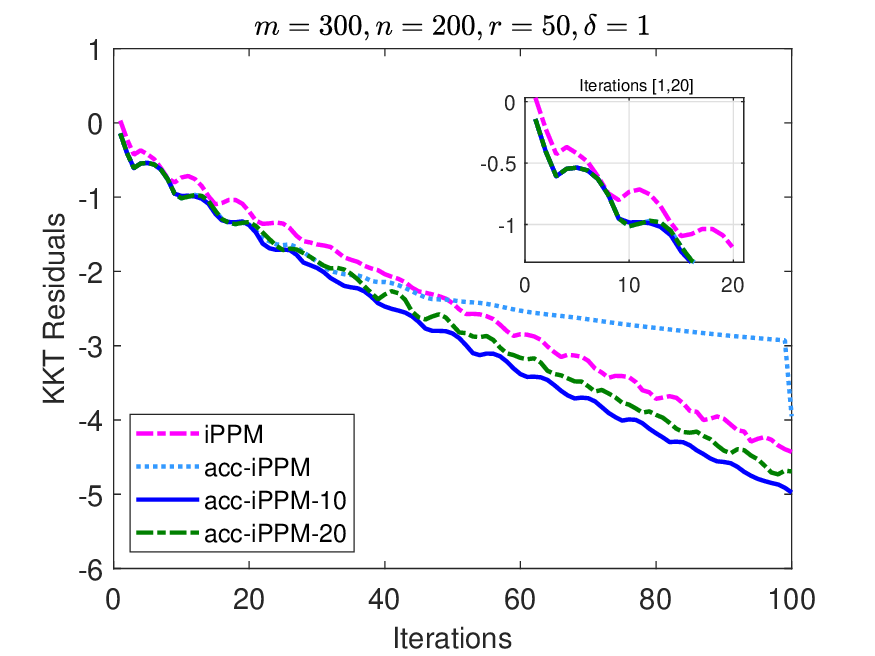}
\hfill
\includegraphics[scale=0.33]{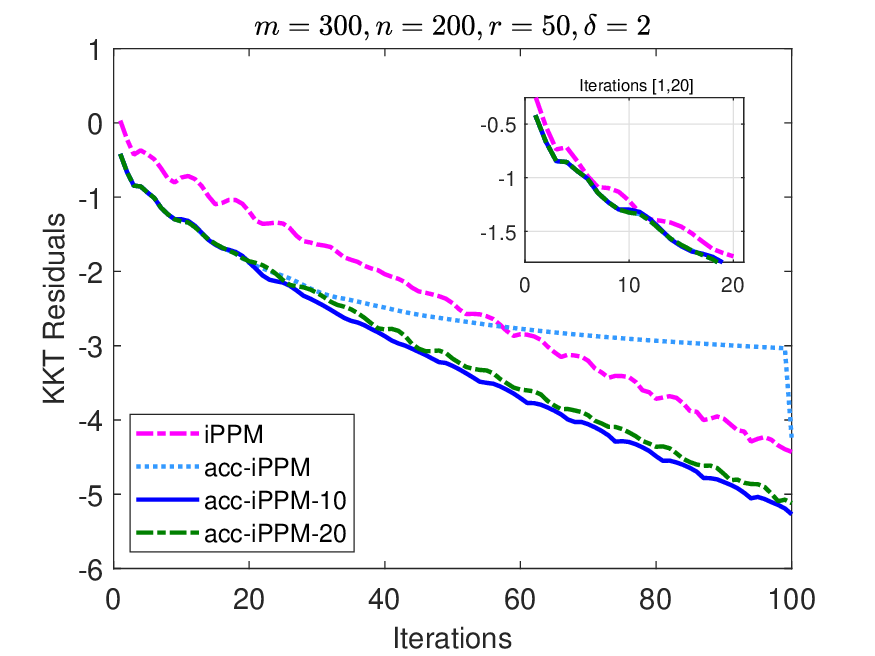}
\hfill
\includegraphics[scale=0.33]{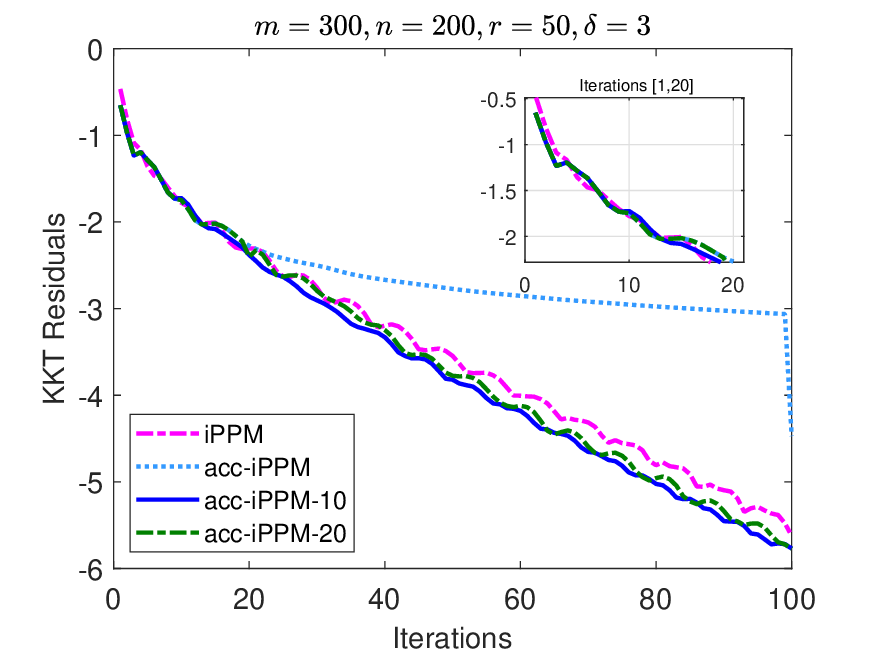}
\caption{Comparison of the KKT residuals for iPPM and HiPPM with and without restarting.}
\label{fig-KKTres}
\end{figure}

We first compare the performance of the iPPM, the HiPPM without restarting, and the HiPPM with restarting every 10 and 20 iterations, respectively. The performance is evaluated in terms of the KKT residual, measured on a base-10 logarithmic scale, defined as
\[
r(z) = \log_{10}\!\left(\|z-P_c(z-A^T(Az-B))\|\right).
\]
The numerical results are summarized in Figure~\ref{fig-KKTres}. To highlight the effect of the Halpern accelerate iteration, we present a magnified view of the first 20 iterations. As shown in the figure, the HiPPM exhibits faster convergence than the iPPM during the early iterations. However, as the number of iterations increases, the convergence speed of the HiPPM without restarting may deteriorate due to the influence of the fixed anchor point. In contrast, the restarting strategy effectively alleviates this issue and enables the HiPPM to consistently achieve a faster convergence rate than the iPPM.

Since Theorem~\ref{Th:3.1} characterizes the convergence behavior of the quantity $\|Z^k-P_c(Z^k)\|$ under the prescribed stopping criterion \eqref{eq:computable_stop}, we first note that
\[
\|z^k-P_c(z^k)\|
\leq \|z^k-\bar{z}^k\|+\|\bar{z}^k-P_c(z^k)\|
\leq \|z^k-\bar{z}^k\|+\varepsilon_k.
\]

\begin{figure}[H]
\centering
\includegraphics[scale=0.33]{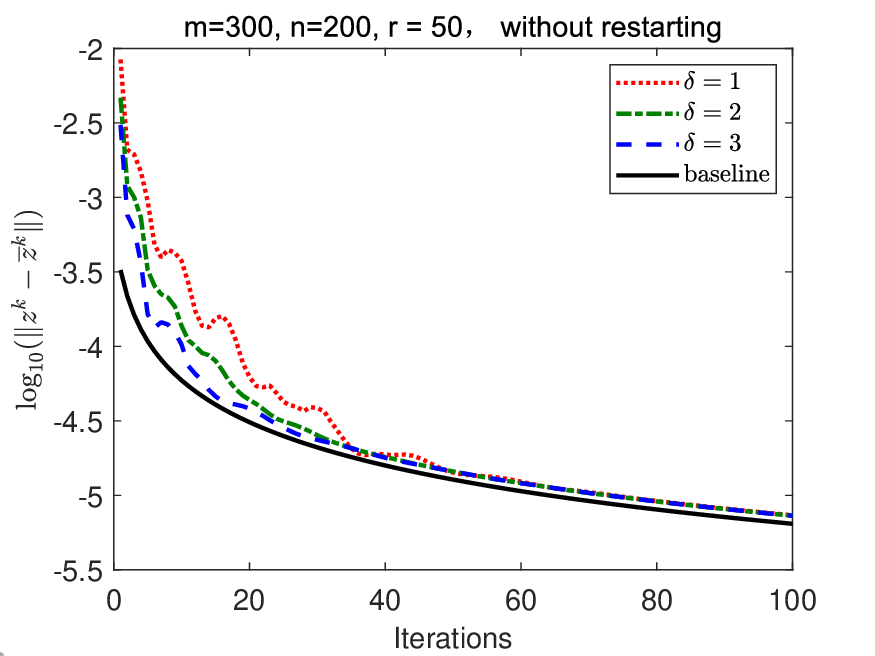}
\includegraphics[scale=0.33]{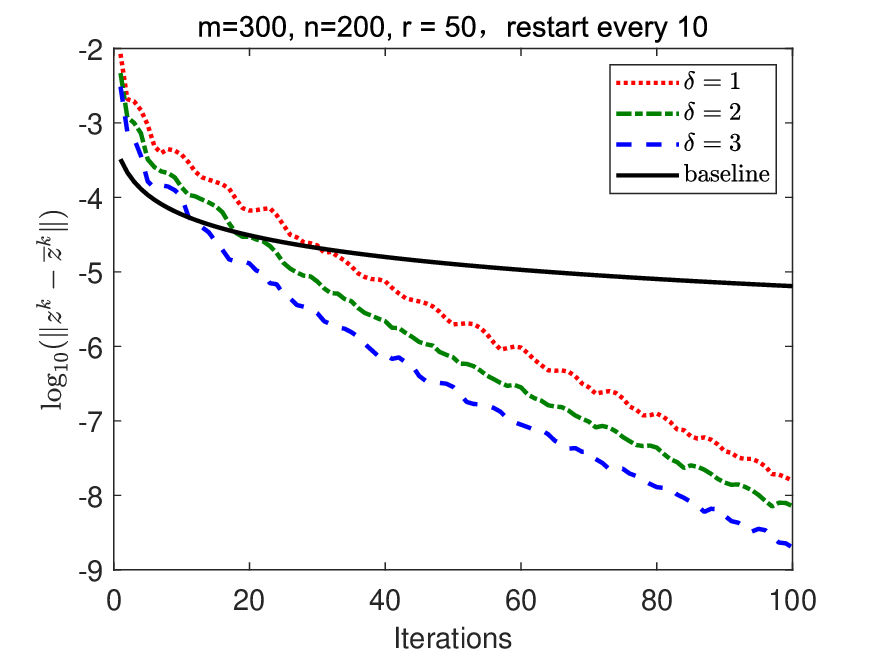}
\includegraphics[scale=0.33]{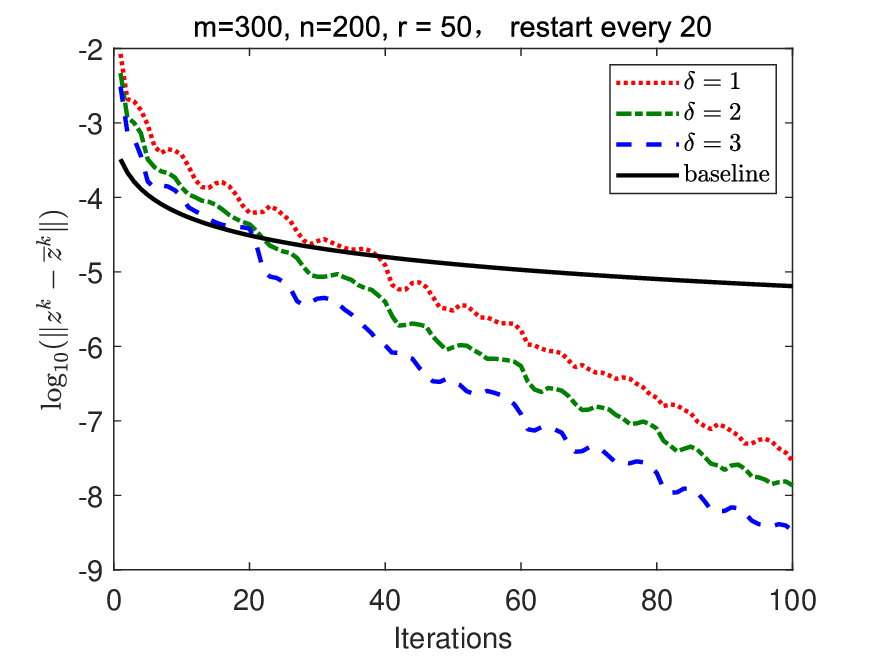}
\caption{Decay of $\|z^k-\overline{z}^k\|$ for different values of $\delta\in\{1,2,3\}$, compared with the baseline $\alpha/(k+1)$.}
\label{fig-normrate}
\end{figure}

Therefore, under the adopted stopping rule, it suffices to investigate the convergence rate of $\|z^k-\bar{z}^k\|$. We compare its behavior for different values of $\delta\in\{1,2,3\}$ against the baseline ${\rm log}_{10}(\alpha/(k+1))$, where $\alpha$ is defined as the average of $\|z^0-z^*\|$ over all runs, and $z^*$ denotes an approximate solution obtained when the stopping condition $r(z)<{-6}$ is satisfied in each run. The corresponding results are reported in Figure~\ref{fig-normrate}. The numerical results show that, in the absence of restarting, the decay rate is initially slower than the baseline, and the case $\delta=3$ is close to the baseline. As the iteration count increases, the convergence behavior for all three values of $\delta$ becomes comparable to that of the baseline. Moreover, the restarting strategy is observed to significantly accelerate convergence.

\subsection{Fast linear convergence rate under regularity conditions}

We consider the following sparse logistic regression problem:
\begin{equation}\label{model:sparse_logistic}
\min_{z\in\mathbb{R}^{n}} \; h(Az) + \lambda\|z\|_1,
\end{equation}
where, for a given vector $b\in\mathbb{R}^m$, the function $h:\mathbb{R}^{m}\to\mathbb{R}$ is defined as
\[
h(y) = \sum_{i=1}^{m} \log\bigl(1+e^{-b_i y_i}\bigr), \qquad \forall\, y\in\mathbb{R}^m.
\]

Let $T(z):=\partial\bigl(h(Az)+\lambda\|z\|_1\bigr)$ denote the associated maximal monotone operator. According to \cite[Theorem~3.3]{li2018highly}, the operator $T$ satisfies an error bound condition at a point $0$ with modulus $a\geq 0$. Specifically, there exists a constant $\bar{\varepsilon}>0$ such that, for any $z\in\mathbb{R}^n$ satisfying ${\rm dist}(0,T(x))\leq \bar{\varepsilon}$, one has
$$
{\rm dist}\bigl(z,T^{-1}(0)\bigr)
\leq a\,{\rm dist}\bigl(0,T(x)\bigr).
$$
Moreover, it follows from \cite[Theorem~3.2]{li2018highly} that, under the above regularity condition, the convergence results established in Theorem~\ref{th4.1} remain valid.

In this experiment, we evaluate the fast linear convergence behavior predicted by Theorem~\ref{th4.1} using two benchmark binary classification datasets from the UCI repository, namely \textit{colon-cancer} and \textit{duke}. The \textit{colon-cancer} dataset consists of 62 samples with 2000 features, while the \textit{duke} dataset contains 44 samples with 7129 features. We set the regularization parameter to $\lambda=10^{-4}$ and incorporate the Halpern iteration into the implementation based on the publicly available code\footnote{\url{https://github.com/linmeixia/exclusive-lasso-solver}}. The subproblems are solved using a semismooth Newton method, and the proximal parameter $c_k$ is chosen to increase progressively along the iterations.

\begin{figure}[H]
\centering
\includegraphics[scale=0.48]{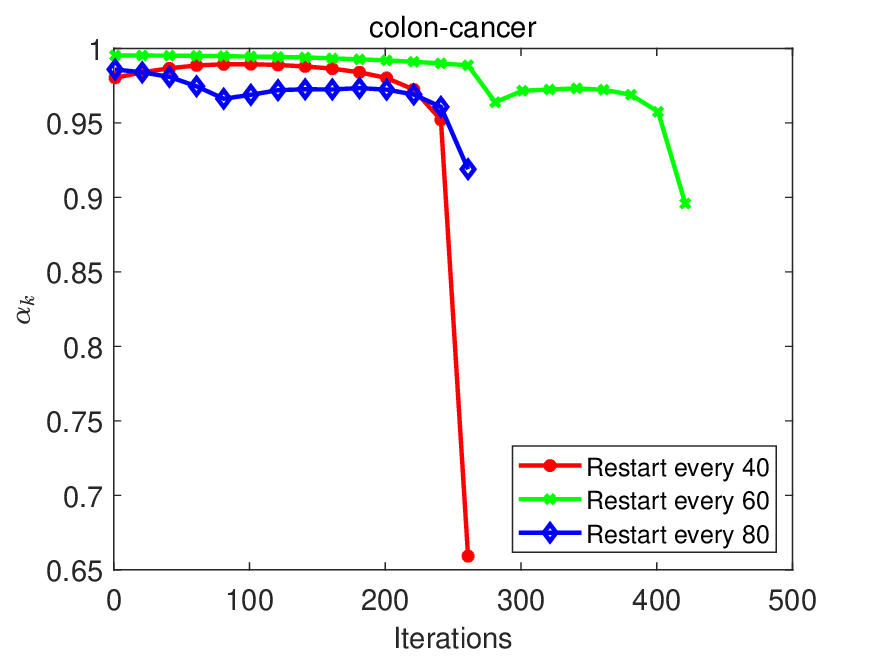}
\hfill
\includegraphics[scale=0.48]{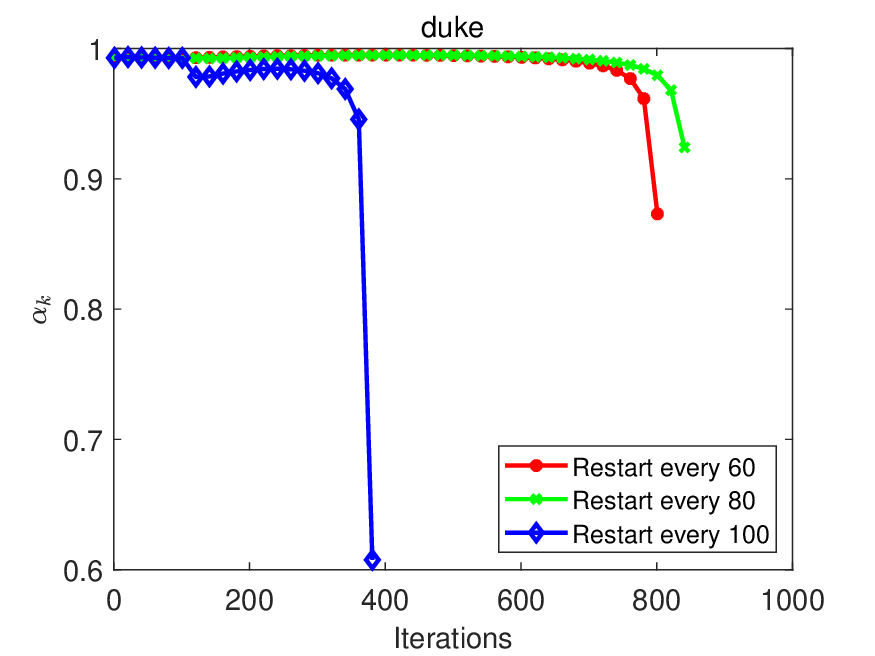}
\caption{Solving problem with regularity condition, $\alpha_k ={\|\bar{z}^k - \bar{z}\|}/{\|z^k - \bar{z}\|}$}
\label{fig-FLrate}
\end{figure}

We report the evolution of the quantity
$
\alpha_k := {\|\bar{z}^k-\bar{z}\|}/{\|z^k-\bar{z}\|}
$
for both datasets under different restarting frequencies in Figure~\ref{fig-FLrate}. The numerical results show that $\alpha_k<1$ for all iterations, which implies the parameter $\vartheta_k$ in \eqref{eq:R3.2} is less than $1$. Furthermore, as the parameter $c_k$ increases, the values of $\alpha_k$ decrease accordingly and thus confirms the fast linear convergence behavior, which is consistent with the theoretical predictions of Theorem~\ref{Th:4.1}.

\section{Conclusion}\label{sec6:conclusion}

In this paper, we studied a Halpern-type acceleration of the inexact proximal point method for maximal monotone inclusion problems in Hilbert spaces. We established its global convergence and developed a unified framework for convergence rate analysis, showing that the squared fixed-point residual converges at an $\mathcal{O}(1/k^2)$ rate under mild inexactness conditions. Moreover, we proved that the proposed Halpern-accelerated scheme achieves fast linear convergence rate under some regularity assumption.  We further extended the analysis to constrained convex optimization via the augmented Lagrangian framework, deriving convergence rate and complexity results for the resulting accelerated inexact augmented Lagrangian method. From a practical perspective, the numerical performance can be further improved by restarting strategies. However, the restart mechanism adopted in this work is heuristic and the design of theoretically justified and adaptive restarting rules remains an important direction for future research.

\end{document}